\documentclass[11pt]{article}
\usepackage{amsfonts,mathrsfs,amssymb,amsthm,mathptm}
\usepackage{amsmath,amscd}
\usepackage{mathptm,pslatex}
\usepackage{multirow}
\usepackage{color}
\usepackage[all]{xy}
\usepackage{url}
\usepackage{cite}
\usepackage{exscale}
\usepackage{relsize}
\usepackage{bbm}
\begin{document}
	\newtheorem{Def}{Definition}[section]
	\newtheorem{Bsp}[Def]{Example}
	\newtheorem{Prop}[Def]{Proposition}
	\newtheorem{Theo}[Def]{Theorem}
	\newtheorem{Lem}[Def]{Lemma}
	\newtheorem{Koro}[Def]{Corollary}
    \newtheorem{Ques}[Def]{Question}
	\theoremstyle{definition}
	\newtheorem{Rem}[Def]{Remark}
	
	\newcommand{\add}{{\rm add}}
	\newcommand{\con}{{\rm con}}
	\newcommand{\gd}{{\rm gl.dim}}
	\newcommand{\dm}{{\rm domdim}}
	\newcommand{\tdim}{{\rm dim}}
	\newcommand{\E}{{\rm E}}
	\newcommand{\Mor}{{\rm Morph}}
	\newcommand{\End}{{\rm End}}
	\newcommand{\ind}{{\rm ind}}
	\newcommand{\rsd}{{\rm res.dim}}
	\newcommand{\rd} {{\rm rep.dim}}
	\newcommand{\ol}{\overline}
	\newcommand{\overpr}{$\hfill\square$}
	\newcommand{\rad}{{\rm rad}}
	\newcommand{\soc}{{\rm soc}}
	\renewcommand{\top}{{\rm top}}
	\newcommand{\stp}{{\mbox{\rm -stp}}}
	\newcommand{\pd}{{\rm projdim}}
	\newcommand{\id}{{\rm injdim}}
	\newcommand{\fld}{{\rm flatdim}}
	\newcommand{\fdd}{{\rm fdomdim}}
	\newcommand{\Fac}{{\rm Fac}}
	\newcommand{\Gen}{{\rm Gen}}
	\newcommand{\fd} {{\rm findim}}
	\newcommand{\Fd} {{\rm Findim}}
	\newcommand{\Pf}[1]{{\mathscr P}^{<\infty}(#1)}
	\newcommand{\DTr}{{\rm DTr}}
	\newcommand{\cpx}[1]{#1^{\bullet}}
	\newcommand{\D}[1]{{\mathscr D}(#1)}
	\newcommand{\Dz}[1]{{\mathscr D}^+(#1)}
	\newcommand{\Df}[1]{{\mathscr D}^-(#1)}
	\newcommand{\Db}[1]{{\mathscr D}^b(#1)}
	\newcommand{\C}[1]{{\mathscr C}(#1)}
	\newcommand{\Cz}[1]{{\mathscr C}^+(#1)}
	\newcommand{\Cf}[1]{{\mathscr C}^-(#1)}
	\newcommand{\Cb}[1]{{\mathscr C}^b(#1)}
	\newcommand{\Dc}[1]{{\mathscr D}^c(#1)}
	\newcommand{\K}[1]{{\mathscr K}(#1)}
	\newcommand{\Kz}[1]{{\mathscr K}^+(#1)}
	\newcommand{\Kf}[1]{{\mathscr  K}^-(#1)}
	\newcommand{\Kb}[1]{{\mathscr K}^b(#1)}
	\newcommand{\modcat}{\ensuremath{\mbox{{\rm -mod}}}}
	\newcommand{\Modcat}{\ensuremath{\mbox{{\rm -Mod}}}}
	
	\newcommand{\stmodcat}[1]{#1\mbox{{\rm -{\underline{mod}}}}}
	\newcommand{\pmodcat}[1]{#1\mbox{{\rm -proj}}}
	\newcommand{\imodcat}[1]{#1\mbox{{\rm -inj}}}
	\newcommand{\Pmodcat}[1]{#1\mbox{{\rm -Proj}}}
	\newcommand{\Imodcat}[1]{#1\mbox{{\rm -Inj}}}
	\newcommand{\PI}[1]{#1\mbox{{\rm -prinj}\,}}
	\newcommand{\opp}{^{\rm op}}
	\newcommand{\otimesL}{\otimes^{\rm\mathbb L}}
	\newcommand{\rHom}{{\rm\mathbb R}{\rm Hom}\,}
	\newcommand{\projdim}{\pd}
	\newcommand{\Hom}{{\rm Hom}}
	\newcommand{\Coker}{{\rm Coker}}
	\newcommand{ \Ker  }{{\rm Ker}}
	\newcommand{ \Cone }{{\rm Con}}
	\newcommand{ \Img  }{{\rm Im}}
	\newcommand{\Ext}{{\rm Ext}}
	\newcommand{\StHom}{{\rm \underline{Hom}}}
	
	\newcommand{\gm}{{\rm _{\Gamma_M}}}
	\newcommand{\gmr}{{\rm _{\Gamma_M^R}}}
	
	\def\vez{\varepsilon}\def\bz{\bigoplus}  \def\sz {\oplus}
	\def\epa{\xrightarrow} \def\inja{\hookrightarrow}
	
	\newcommand{\lra}{\longrightarrow}
	\newcommand{\llra}{\longleftarrow}
	\newcommand{\lraf}[1]{\stackrel{#1}{\lra}}
	\newcommand{\llaf}[1]{\stackrel{#1}{\llra}}
	\newcommand{\ra}{\rightarrow}
	\newcommand{\dk}{{\rm dim_{_{k}}}}
	
	\newcommand{\colim}{{\rm colim\, }}
	\newcommand{\limt}{{\rm lim\, }}
	\newcommand{\Add}{{\rm Add }}
	\newcommand{\Tor}{{\rm Tor}}
	\newcommand{\Cogen}{{\rm Cogen}}
	\newcommand{\Tria}{{\rm Tria}}
	\newcommand{\tria}{{\rm tria}}
	
{\Large \bf
\begin{center}
Structure of centralizer algebras
\end{center}}
	
\medskip
\centerline{\textbf{Changchang Xi$^*$ and Jinbi Zhang}}
	
\renewcommand{\thefootnote}{\alph{footnote}}
\setcounter{footnote}{-1} \footnote{$^*$Corresponding author's
		Email: xicc@cnu.edu.cn; Fax: 0086 10 68903637.}
\renewcommand{\thefootnote}{\alph{footnote}}
\setcounter{footnote}{-1} \footnote{2010 Mathematics Subject
		Classification: Primary 16S50, 15B33, 16U70, 15A27; Secondary
		 20C05, 16W22, 11C20.}
\renewcommand{\thefootnote}{\alph{footnote}}
\setcounter{footnote}{-1} \footnote{Keywords: Centralizer algebra; Cellular algebra; Frobenius extension; Jordan block; Matrix algebra; Symmetric group. }
	
\begin{abstract}
Given an $n\times n$ matrix $c$ over a unitary ring $R$, the centralizer of $c$ in the full $n\times n$ matrix ring $M_n(R)$ is called a principal centralizer matrix ring, denoted by $S_n(c,R)$. We investigate its structure and prove: $(1)$ If $c$ is an invertible matrix with a $c$-free point, or if $R$ has no zero-divisors and $c$ is a Jordan-similar matrix with all eigenvalues in the center of $R$, then $M_n(R)$ is a separable Frobenius extension of $S_{n}(c,R)$ in the sense of Kasch. $(2)$ If $R$ is an integral domain and $c$ is a Jordan-similar matrix, then $S_n(c,R)$ is a cellular $R$-algebra in the sense of Graham and Lehrer. In particular, if $R$ is an algebraically closed field and $c$ is an arbitrary matrix in $M_n(R)$, then $S_n(c,R)$ is always a cellular algebra, and the extension $S_n(c,R)\subseteq M_n(R)$ is always a separable Frobenius extension.
	\end{abstract}

	
\section{Introduction}
Let $R$ be a unitary (associative) ring and $C$ a nonempty subset of $R$. The \emph{centralizer} of $C$ in $R$ is a subring of $R$ defined by $$S(C,R): =\{r\in R \mid cr=rc \;\mbox{ for \;all\;} c\in C\}.$$ The center of $R$ is $S(R,R)$, denoted by $Z(R)$.  If $C=\{c\}$ is a singleton set, then $S(c,R):=S(\{c\},R)$ is called a \emph{principal centralizer} ring. Clearly, $S(C,R)=\bigcap\limits_{c\in C}S(c,R)$ for any nonempty subset $C$ of $R$. This means that, to understand the entire ring $S(C,R)$, it may be useful to study first the individual ring $S(c,R)$ for each $c\in C$.

The centralizer rings $S(C,R)$ have been related to invariant rings. Let $G$ be a subgroup of the group Aut$(R)$ of automorphisms of the ring $R$. For any $g\in$ Aut$(R)$, the action of $g$ on $R$ is denoted by $r\mapsto r^g$ for $r\in R$. An $r\in R$ is said to be \emph{$G$-free} if $r^g\neq r$ for $g\ne id$, that is, the stabilizer group st$_G(r):=\{g\in G \mid r^g=r\}$ of $r$ under $G$ is trivial. An element $g\in$ Aut$(R)$ is said to be \emph{inner} if there exists a unit $s_g\in R$ such that $r^g=s_g^{-1}rs_g$ for all $r\in R$. The fixed ring (or invariant ring) of $R$ under $G$ is $R^G :=\{r\in R \mid r^g=r\; \mbox{ for\; all } g\in G\}$. If all elements of $G$ are inner and $C=\{s_g\in R\mid g\in G\}$, then $S(C,R)=R^G$. The study of fixed rings has a long history (see \cite{W}). They are investigated by many authors, we refer to \cite{Montgomery1980} and the references therein for more information on fixed rings with $G$ finite groups.

For a positive integer $n$, $M_n(R)$ denotes the full matrix ring of all $n\times n$ matrices over $R$, and  $GL_n(R)$ the general linear group of all invertible $n\times n$ matrices over $R$.  For a nonempty subset $C$ of $M_n(R)$, the ring $S_n(C,R): = S(C,M_n(R))$ is called  a \emph{centralizer matrix ring} over $R$ of degree $n$. For a matrix $c\in M_n(R)$, the ring $S_n(c,R): = S(c,M_n(R))$ is called  a \emph{principal centralizer matrix ring} over $R$ of degree $n$, and the extension $S_n(c,R)\subseteq M_n(R)$ is called a {\em principal centralizer matrix extension}. They are the objectives of this note. Typical examples of principal centralizer matrix rings and extensions include centrosymmetric matrix algebras (see \cite{Weaver, Xi2020}).

Recall that a matrix in $ M_n(R)$ is a \emph{Jordan block} if it is of the form
$$\begin{pmatrix}  r & 1         & \cdots & 0& 0\\
              0 & r         &  \ddots & 0 &0 \\
         \vdots &\vdots     & \ddots & \ddots & \vdots\\
          0     & 0          & \cdots & r & 1 \\
          0   &   0& \cdots  & 0 & r\\
\end{pmatrix}_{n\times n}$$
where $r\in R$ is called the \emph{eigenvalue} of the matrix. A matrix $a=(a_{ij})\in M_n(R)$ is called a \emph{Jordan-block matrix} if it is a diagonal-block matrix with each block in the main diagonal being a Jordan block. In this case, we may suppose that $a$ has $t$ distinct eigenvalues in $R$, say $r_1, \cdots, r_t$, and that, for each eigenvalue $r_i$, there are $s_i$ Jordan-blocks $J_{ij}$ of distinct sizes $\lambda_{ij}$ with the eigenvalue $r_i$, $1\le j\le s_i$, such that $J_{ij}$ appears $b_{ij}$ time and $\lambda_{i1}>\lambda_{i2}>\cdots> \lambda_{is_i}$ for $1\le i\le t$. The set $\{(\lambda_{11}, \lambda_{12},\cdots,\lambda_{1s_1}),(\lambda_{21},\lambda_{22},\cdots, \lambda_{2s_2}),\cdots, (\lambda_{t1},\lambda_{t2},\cdots,\lambda_{ts_t})\}$ is called the \emph{block type} of $a$.
A matrix $c\in M_n(R)$ is called a \emph{Jordan-similar matrix} if it is similar to a Jordan-block matrix $a$ by a matrix in $GL_n(R)$, that is, there is a matrix $u\in GL_n(R)$ such that $a=ucu^{-1}$ is a Jordan-block matrix in $M_n(R)$. In this case, the \emph{block type} of $c$ is defined to be the block type of $a$.
If $R$ is an algebraically closed field, then every square matrix over $R$ is Jordan-similar (for example, see \cite[\uppercase\expandafter{\romannumeral7}.7]{G}).

In this note, we investigate cellular structures in the sense of Graham and Lehrer (see \cite{Graham1996}), and Frobenius extension properties in the sense of Kasch (see \cite{Kasch1961}), of centralizer matrix rings. This is approached by combining methods in matrix theory with the ones in combinatorics and representation theory of algebras. The idea to combine them in proofs seems to be new.

Our first main result points out a cellular structure of principal centralizer matrix algebras.

\begin{Theo} \label{MAIN2} Let $R$ be an integral domain and $c$ a Jordan-similar matrix in $M_n(R)$. Then

$(1)$ $S_n(c,R)$ is a cellular $R$-algebra.		
				
$(2)$ Suppose that $R$ is a field and $c$ is of the block type $\{(\lambda_{11}, \lambda_{12},\cdots,\lambda_{1s_1}),\cdots, (\lambda_{t1},\lambda_{t2},\cdots,\lambda_{ts_t})\}.$ Then $S_n(c,R)$ is a quasi-hereditary algebra if and only if $\lambda_{i1}=s_i$ for $1\le i \le t$.
\end{Theo}

Next, we establish Frobenius extensions of centralizer matrix algebras included in full matrix algebras.

\begin{Theo}\label{MAIN1} Let $R$ be a unitary ring.

$(1)$ If $G$ is a finite subgroup of $GL_n(R)$ with a $G$-free point, then $S_{n}(G,R)\subseteq M_n(R)$ is a separable Frobenius extension. Moreover, if $|G|$ is additionally invertible in $R$, then the extension is also split.

$(2)$ Suppose that $R$ has no zero-divisors and $c\in M_n(R)$ is a Jordan-similar matrix with all eigenvalues in $Z(R)$. Then $S_n(c,R)\subseteq M_n(R)$ is a separable Frobenius extension. Moreover, the extension $S_n(c,R)\subseteq M_n(R)$ is split if and only if $c$ is similar to a matrix of the form ${\rm diag}(r_1I_{n_1},r_2I_{n_2},\cdots,r_tI_{n_t})$ with $\sum_in_i=n$, $r_i\in Z(R)$ and $r_i\neq r_j$ for $1\le i,j\le t$.
\end{Theo}

Thus, if $R$ is an algebraically closed field, then every principal centralizer matrix extension over $R$ is a Frobenius extension, and every principal centralizer matrix $R$-algebra is cellular. Moreover, the number of non-isomorphic simple modules of such a cellular algebra can be described combinatorially by the data of Jordan forms (see Corollary \ref{quasi-hered}). It is surprising that, in general, if $C$ has more than one element or $G$ does not have any free point, then Theorem \ref{MAIN1} is no longer true. This is demonstrated by examples in Section \ref{FROB}.

Consequently, we have the corollary.
	
\begin{Koro}\label{MAINCOR}
Let $G$ be a finite group, $k$ an algebraically closed field such that its characteristic does not divide the order of $G$, and $c$ an element in the group algebra $kG$ of $G$ over $k$. Then $S(c, kG)$ is a cellular algebra and $S(c, kG) \subseteq kG$ is a separable Frobenius extension.
\end{Koro}
	
Theorem \ref{MAIN1}(1) generalizes \cite[Theorem 3.1(3)]{Xi2020} substantially, while Corollary \ref{MAINCOR} extends \cite[Theorem 3.3]{Xi2020} in case of $R$ being an algebraically closed field. As principal centralizer matrix algebras are much more general and complicated than centrosymmetric matrix algebras, our approach in this note is different from the one in \cite{Xi2020}. For example, the involution used for the cellular structure of principal centralizer matrix algebras is completely different from the matrix transpose used in \cite{Xi2020}.

The paper is organized as follows: In Section \ref{CMA}, we fix notation and develop basic facts on  centralizer matrix algebras. In Section \ref{FROB}, we first recall the definition of and some known results on Frobenius extensions, and then prove Theorem \ref{MAIN1}. In Section \ref{CELL}, we show Theorem \ref{MAIN2} and Corollary \ref{MAINCOR} after recalling the notion of cellular algebras. The section ends with a few open questions.
	
\section{Principal centralizer matrix algebras\label{CMA}}
In this section we discuss basic properties of centralizer matrix algebras.

Let $m,n$ be positive integers. We write  $[n]$ for the set $\{1,2,\dots,n\}$.

Let $R$ be a unitary ring ( that is, an associative ring with identity), we denote by $\rad(R)$ the Jacobson radical of $R$, by $M_{m\times n}(R)$ the set of all $m\times n$ matrices over $R$ and by $e_{ij}$ the matrix units of $M_{m\times n}(R)$, with $i\in[m], j\in [n]$. We write
$M_n(R)$ for $M_{n\times n}(R)$ and $I_n$ for the identity matrix in $M_n(R)$. For a matrix $a\in M_{m\times n}(R)$, we denote by $a'$ the transpose of $a$.

By an $R$-module we means a left $R$-module. For an $R$-module $M$, $\End_R(M)$ stands for the endomorphism ring of $M$.
If $f: X\ra Y$ and $g:Y\ra Z$ are homomorphisms of modules, the composite of $f$ and $g$ will be denoted by $fg$. This means that the image of $x\in X$ under $f$ is written as $(x)f$, instead of $f(x)$.

\begin{Lem}\label{PROPCMA} Let $C$ be a subset of $R$.
		
$(1)$ If $x$ is invertible in $M_n(R)$, then there are isomorphisms of rings: $S_n(xCx^{-1},R)\simeq S_n(C,R)$ and $S_n(Cx,R)\simeq S_n(xC,R)$, where $xCx^{-1}:=\{xcx^{-1}\mid c\in C\}$, $xC:=\{xc\mid c\in C\}$ and $Cx:=\{cx\mid c\in C\}$.
		
$(2)$ If each $c\in C$ is invertible in $M_n(R)$ and $C^{-1}:=\{c^{-1}\mid c\in C\}$, then $S_n(C,R)\simeq S_n(C^{-1},R)$.
		
$(3)$ If $C=\{c^i\in M_n(R)\mid i\in\mathbb{N}\}$, then $S_n(C,R)=S_n(c,R)$.
\end{Lem}
	
{\it Proof.} (2) and (3) are trivial. For (1), the correspondence $r\mapsto x^{-1}rx$ gives rise to an isomorphism of rings not only from $S_n(xCx^{-1},R)$ to $S_n(C,R)$, but also from $S_n(x C,R)$ to $S_n(Cx,R)$. $\square$

\medskip	
Recall that an $n\times n$ matrix $a$ is called \emph{semicirculant} if it has the form
$$
	\begin{pmatrix}
		a_{1}& a_{2}& a_{3} &\dots&a_{n}\\
		0&a_{1}&a_{2}&\dots&a_{n-1}\\
		0&0&a_{1}&\dots&a_{n-2}\\
		\vdots&\vdots&\vdots&\ddots&\vdots\\
		0&0 &0 &\dots&a_{1}
	\end{pmatrix}\in M_n(R).
$$
In this case, we write $a =[a_{1},a_{2},\dots,a_{n}]$. If $J_{n,0}=[0,1,0,\ldots,0]$ stands for the Jordan block of size $n$ with the eigenvalue $0$, then the semicirculant matrix $a=[a_{1},a_{2},\ldots,a_{n}]$ can be represented as a polynomial in $J_{n,0}$:
\begin{eqnarray*}		a=a_{1}J_{n,0}^0+a_{2}J_{n,0}+\cdots+a_{n}J_{n,0}^{n-1}=\sum_{p=1}^{n}a_{n-p+1}
\sum_{j=1}^{p}e_{p-j+1,n-j+1}
=\sum_{p=1}^{n}\sum_{j=1}^{p}a_{n-p+1}e_{p-j+1,n-j+1,}
\end{eqnarray*}
where $J_{n,0}^0$ is understood as the $n\times n$ identity matrix $I_n$.

For $1\le p\le \mbox{min}\{m,n\}$, we define
$$G^p :=\sum_{j=1}^{p}e_{p-j+1,n-j+1}\in M_{m\times n}(R). $$
If $m=n$, then $[a_1,\cdots,a_n]=\sum_{i=1}^na_{n-i+1}G^i$.	

\begin{Lem}\label{lemsemicir} Let $I=[r,1,0,\dots,0]\in M_m(R)$ and $J=[r',1,0,\dots,0]\in M_n(R)$ be Jordan blocks with $r,r'\in Z(R)$. Define $\Lambda :=\{a\in M_{m\times n}(R)\mid Ia=aJ\}$.
		
$(1)$ If $r\neq r'$ and $R$ has no zero-divisors, then $\Lambda = 0$.

$(2)$ If $r=r'$, then $\Lambda$ is a free $R$-module with an $R$-basis $\big\{G^p\mid 1\le p\le \emph{min}\{m,n\}\big\}.$
	\end{Lem}
	
	{\it Proof.} Let $a=(a_{ij})\in M_{m\times n}(R)$ and $a_{i0} = a_{m+1,j}=0$ for $1\leq i\leq m,\;1\leq j\leq n$. Then $Ia=(r a_{ij}+a_{i+1,j})_{m\times n}.$
Similarly, $aJ=(a_{i,j-1}+a_{ij}r')_{m\times n}$.
Due to $r,r'\in Z(R)$, we know that
$Ia=aJ$ if and only if
$$(*)\quad (r-r')a_{ij}= a_{i,j-1}-a_{i+1,j} \mbox{  for all } 1\le i\le m, 1\le j\le n.$$
Thus $a=(a_{ij})\in M_{m\times n}(R)$ lies in $\Lambda$ if and only if $(*)$ holds.

(1) Assume $r\ne r'$. Then $(r-r')a_{m1}= a_{m 0}- a_{m+1,1}=0$. Since $R$ has no zero-divisors, we have $a_{m1}=0$. It follows from $(r-r')a_{i1}= a_{i 0}-a_{i+1,1}=-a_{i+1,1}$ that $a_{i 1}=0$ for $1\le i\le m$. Similarly, $(r-r')a_{i 2}= a_{i 1}-a_{i+1,2}=-a_{i+1,2}$ for $1\le i\le m$. This implies  $a_{m 2}=a_{m-1,2}= \cdots = a_{12}=0$. Continuing this argument, we get $a_{ij}=0$ for $1\leq i\leq m$, $3\leq j\leq n$. Hence $a=0$ and $\Lambda=\{0_{m n}\}$.
	
(2) Assume $r=r'$. Then $a=(a_{ij})\in M_{m\times n}(R)$ lies in $\Lambda$ if and only if $a_{i,j-1}=a_{i+1,j}\mbox{ for }1\le i\le m, 1\le j\le n.$
Thus $a_{ij}=a_{i-1,j-1}=\cdots=a_{i-j,0}=0$ for $1\le i-j\le m-1$ and $a_{ij}=a_{i+1,j+1}=\cdots=a_{m+1,j+m+1-i}=0$ for $m-n+1\le i-j\le m-1$. That is, $a_{ij}=0$ for $\mbox{min}\{1,m-n+1\}\le i-j\le m-1 $.

Let $l:=\mbox{min}\{m,n\}.$ For $1-n\le i-j\le \mbox{min}\{0,m-n\}$, it follows from ($*$) that $a_{ij}=a_{i-1,j-1}=\cdots=a_{1,j-i+1}$. Let $p:=i-j+n$ and $u:=n-j+1$. Then $1\le p\le l$, $i=p-u+1$, $j=n-u+1$ and $1\le u\le n$. It follows from $0\le p-u=i-1\le m-1$ that $a_{1,n-p+1}=a_{p-u+1,n-u+1}$ for $1\le u\le p\le l$. Thus $a=\sum_{p=1}^{l}a_{1,n-p+1}\sum_{u=1}^{p}e_{p-u+1,n-u+1}=\sum_{p=1}^{l}a_{1,n-p+1}G^p$. Hence  $a\in \Lambda$ can be written as an $R$-linear combination of $\{G^p\mid 1\le p\le l\}$.

We show $G^p\in \Lambda$ for $1\le p\le l$, that is, $(G^p)_{i,j-1}=(G^p)_{i+1,j}$ for $1\le i\le m$ and $1\le j\le n$. In fact,
\begin{align*}
  e_{i+1,i+1}G^pe_{jj} & =e_{i+1,i+1}\sum_{v=1}^{p}e_{p-v+1,n-v+1}e_{jj}
    =\sum_{v=1}^{p}\delta_{i+1,p-v+1}\delta_{n-v+1,j}e_{i+1,j,} \\
  e_{ii}G^pe_{j-1,j-1} & =e_{ii}\sum_{u=1}^{p}e_{p-u+1,n-u+1}e_{j-1,j-1}
    =\sum_{u=1}^{p}\delta_{i,p-u+1}\delta_{n-u+1,j-1}e_{i,j-1} \\
   & =\sum_{u=1}^{p}\delta_{i+1,p-u+2}\delta_{n-u+2,j}e_{i,j-1}
    =\sum_{v=0}^{p-1}\delta_{i+1,p-v+1}\delta_{n-v+1,j}e_{i,j-1}\\ &=\sum_{v=1}^{p}\delta_{i+1,p-v+1}\delta_{n-v+1,j}e_{i,j-1.}
\end{align*}
Thus $(G^p)_{i,j-1}=\sum_{v=1}^{p}\delta_{i+1,p-v+1}\delta_{n-v+1,j}=(G^p)_{i+1,j.}$ This implies that $G^p\in\Lambda$ and $\{G^p\mid 1\le p\le l\}$ is an $R$-generating set of $\Lambda.$

Moreover, $\{G^p\mid 1\le p\le l\}$ is an  $R$-linear independent set.
Indeed, $(p-u+1,n-u+1)=(p'-u'+1,n-u'+1)$ if and only if $(p,u)=(p',u')$.
Since $1 \le p-u+1 \le l\le m$ and $1\le n-u+1\le n$, there holds the inclusion $\{e_{p-u+1,n-u+1} \mid 1 \le p \le l, 1\le u\le p\}\subseteq \{e_{pu}\mid1\le p\le m,1\le u\le n\}$.
As the matrix units $\{e_{pu}\mid1\le p\le m,1\le u\le n\}$ is an $R$-basis of $M_{m\times n}(R)$, we know that if $\sum_{p=1}^{l}a_pG^p=\sum_{p=1}^{l}\sum_{u=1}^{p}a_p e_{p-u+1,n-u+1}=0$ for $a_p\in R$, then $a_p=0$ for $1\le p\le l$, and therefore $\{G^p\mid 1\le p\le l\}$ is a set of $R$-linear independent elements.
Hence $\{G^p\mid 1\le p\le l\}$ is an $R$-basis of $\Lambda$. $\square$

\medskip
In Lemma \ref{lemsemicir}(2), the basis elements $G^p$ does not involve the value $r$. Hence $\Lambda$ is independent of the choice of $r$.
In fact, if we write $I=rI_m+J_{m,0}$ and $J=rI_n+J_{n,0}$, then $Ia=aJ$ for all $a\in M_{m\times n}(R)$ is equivalent to saying that $J_{m,0}a=aJ_{n,0}$ for all $a\in M_{m\times n}(R)$.
Thus $\Lambda$ is independent of the choice of $r$.

A special case of Lemma \ref{lemsemicir} is $c:=I=J=[r,1,0,\dots,0]\in M_n(R)$ for some $r\in Z(R)$. Then $S_n(c,R)$ is the set of all semicirculant matrices in $M_n(R)$. Clearly, $S_n(c,R)\simeq R[x]/(x^n)$ as rings. Moreover, if $R$ is a local ring, that is, the set of all non-units in $R$ is an ideal of $R$, then $S_n(c,R)$ is a local ring.

\begin{Bsp}\label{ex1}{\rm
Let $c_1=[1,1,0]\in M_3(R)$, $c_2=[1,1]\in M_2(R)$ and $c=$diag$(c_1,c_2)\in M_5(R)$. Then $S_5(c,R)$ is a free $R$-module of rank $9$ by Lemma \ref{lemsemicir}(2). If $R$ is a field, then $S_5(c,R)$ is isomorphic to an algebra given by the quiver with relations:
$$
\xymatrix{
1\bullet\ar@<2.5pt>[r]^{\alpha} &\bullet 2,\ar@<2.5pt>[l]^{\beta}
  &  \beta\alpha\beta\alpha=0.}
$$ Here the vertices $1$ and $2$ correspond to the primitive idempotents $f_1=e_{11}+e_{22}+e_{33}$ and $f_2=e_{44}+e_{55}$ in $S_5(c,R)$, respectively, and the compose $\alpha\beta$ of two arrows $\alpha$ and $\beta$ means that $\alpha$ comes first and then $\beta$ follows, that is, $\alpha\beta$ is a path of length $2$ from the vertex $1$ to itself. Thus $S_5(c,R)$ is a cellular algebra with Cartan matrix $\begin{pmatrix} 3&2\\2 &2\end{pmatrix}$ (see \cite[Theorem 6.1]{XX}). Due to $S_3(c_1,R)\simeq R[X]/(X^3)$ and $S_2(c_2,R)\simeq R[X]/(X^2)$, we know $S_5(c,R)\not\simeq S_3(c_1,R)\times S_2(c_2,R).$
This example shows that the study of $S_n(c,R)$ related to $c$ cannot be reduced to the one related to each of Jordan blocks $c_i$. Generally, the centralizer of $c$ in $ M_{m+n}(R)$ does not coincide with the centralizer of $c$ in diag$\big(M_m(R), M_n(R)\big)$.
}\end{Bsp}

Now, we assume that $c$ is a \textbf{Jordan-block matrix with the same eigenvalues in $Z(R)$}.
More precisely, suppose that $c$ has $b_i$ Jordan-blocks $J_i$ of size $\lambda_i$ for $1\le i\le s$ with $\lambda_1 \textgreater \lambda_2 \textgreater \cdots \textgreater \lambda_s$, that is,
$$(\dag)\quad c={\rm diag}(J_{1}^{b_1},J_{2}^{b_2},\dots,J_{s}^{b_s})\in M_{n}(R)$$
with $J_{j}=[r,1,0,\cdots,0]\in M_{\lambda_{j}}(R)$ appearing $b_j$ times for $1\leq j\leq s$ and $r\in Z(R)$.

We define
$$m_0:=0,\;m_i:=\sum_{p=1}^{i}b_p,\; n_{ij}:=\sum_{p=1}^{i-1}b_p\lambda_p+j\lambda_i,\;1\le i\le s,1\le j\le b_i.$$
Then $m_0<m_1<\cdots<m_s$ and $m_s$ is the number of Jordan blocks of $c$ and $n_{sb_{s}}=n$ is the size of the matrix $c$.
For each $i \in [m_s]$, let $g(i)$ be the smallest $g(i) \in [s]$ such that $i\le m_{g(i)}$,  and let  $h(i):=i-m_{g(i)-1}\in [b_{g(i)}]$ and $\theta_{ij}:=\min\{\lambda_{g(i)},\lambda_{g(j)}\} \mbox{ for } j\in [m_s].$ Note that $[n]=\{n_{g(j)h(j)}-\lambda_{g(j)}+w \in \mathbb{N} \mid  1\le j\le m_s,1\le w\le \lambda_{g(j)} \}.$

For each $i \in [m_s]$, we define
$$f_i:=\sum_{p=n_{g(i)h(i)}-\lambda_{g(i)}+1}^{n_{g(i)h(i)}}e_{pp},$$
that is, $f_i$ is the identity matrix corresponding to the $i$-th block in the identity matrix $I_n$. Here  $I_n$ is regarded as a diagonal $m_s\times m_s$ block matrix.
Then $1=\sum_{i=1}^{m_s}f_i$ is a decomposition of $1$ into pairwise orthogonal idempotents of $S_n(c,R)$. Consequently, $S_n(c,R)$ has the matrix decomposition
$$ S_n(c,R)= \begin{pmatrix}
		{f}_1S_n(c,R)f_1 & f_1S_n(c,R)f_2 & \cdots & f_1S_n(c,R) f_{m_s}\\
		f_2S_n(c,R)f_1& f_2 S_n(c,R)f_2 & \cdots &f_2S_n(c,R) f_{m_s}\\
		\vdots & \vdots &\ddots & \vdots\\
		f_{m_s}S_n(c,R)f_1& f_{m_s} S_n(c,R)f_2& \cdots & f_{m_s} S_n(c,R) f_{m_s}
	\end{pmatrix}_{m_s\times m_s}$$
It is easy to see that an $m_s\times m_s$ block matrix $a=(A_{ij})$ with $A_{ij}\in M_{\lambda_{g(i)}\times \lambda_{g(j)}}(R)$ lies in $S_n(c,R)$ if and only if each block $A_{ij}$ satisfies $J_{g(i)}A_{ij}=A_{ij}J_{g(j)}$ for $1\leq i,j \leq m_s$. More precisely, $a\in S_n(c,R)$ is written as a block-matrix form $$a=(A_{ij})=\begin{pmatrix}
		A_{1 1} & A_{1 2} & \dots & A_{1 m_s} \\
		A_{2 1} & A_{2 2} & \dots & A_{2 m_s} \\
		\vdots & \vdots & \ddots &  \vdots \\
		A_{m_s 1} &A_{m_s 2}  & \dots & A_{m_s  m_s}
	\end{pmatrix}_{m_s\times m_s}$$ where the diagonal entry $A_{ii}$ is a $\lambda_{g(i)}\times \lambda_{g(i)}$ semicirculant matrix and the off-diagonal entry $A_{ij}$ is a $\lambda_{g(i)}\times \lambda_{g(j)}$ matrix over $R$, having the block form:$$A_{ij}=\left\{
	\begin{array}{rcl}
		\begin{pmatrix}
			[a_{1},a_{2},\cdots,a_{\lambda_{g(j)}}]\\
			\text{\Large 0}_{\lambda_{g(i)}-\lambda_{g(j)}, \lambda_{g(j)}}
		\end{pmatrix}\quad\qquad&   &if\; \lambda_{g(i)} \textgreater \lambda_{g(j)},\\
		\begin{pmatrix}
			\text{\Large 0}_{\lambda_{g(i)},\lambda_{g(j)}-\lambda_{g(i)}}&[a_{1},a_{2},\cdots,a_{\lambda_{g(i)}}]
		\end{pmatrix}\,&   & if\; \lambda_{g(i)}\textless \lambda_{g(j)},
	\end{array} \right.$$
with all $a_p\in R$. Visually, $$A_{ij}=\left(
	\begin{array}{c|c}
		\text{\Large 0}_{\theta_{ij},\lambda_{g(j)}-\theta_{ij}} & [a_{1},a_{2},\cdots,a_{\theta_{ij}}] \\
		\hline
		\text{\Large 0}_{\lambda_{g(i)}-\theta_{ij},\lambda_{g(j)}-\theta_{ij}} & \text{\Large 0}_{\lambda_{g(i)}-\theta_{ij},\, \theta_{ij}}
	\end{array}
	\right)_.$$

For simplicity,  we set $\Lambda := S_n(c,R)$, $\Lambda_{i j}:=f_iS_n(c,R)f_j$ and $\tilde{\Lambda}_{ij} := \{a\in M_{\lambda_{g(i)}\times \lambda_{g(j)}}(R)\mid J_{g(i)}a=aJ_{g(j)}\}$.
Given $1\le i,j \le m_s$ and $1\le p\le \theta_{ij}$, define
$$G_{ij}^{p}:=\sum_{u=1}^{p}e_{p-u+1,\lambda_{g(j)}-u+1}\in M_{\lambda_{g(i)}\times \lambda_{g(j)}}(R).$$
By Lemma \ref{lemsemicir}(2), $\{G_{ij}^{p}\mid 1\le p\le \theta_{ij}\}$ is an $R$-basis of $\tilde{\Lambda}_{ij}$, which has the property.
\begin{Lem}\label{multbasis1}
$(1)$ If $1\le i,j,k\le m_s,1\le p\le \theta_{ik},1\le q \le \theta_{kj}$, then
$$G_{ik}^pG_{kj}^q=\left\{
	\begin{array}{ll}
		0   & if\;p+q-\lambda_{g(k)}\textless 1, \\ G_{ij}^{p+q-\lambda_{g(k)}}  & if\;p+q-\lambda_{g(k)}\ge 1.
	\end{array} \right.$$
$(2)$ For $1\le i,j,k\le m_s,1\le m\le \theta_{ik},1\le l \le \theta_{kj}$, $x_p,y_q\in R$ for $1\le p\le m$ and $1\le q\le l$, let $$X_{ik}:=\sum_{p=1}^{m}x_{m-p+1}G_{ik}^{p} \in \tilde{\Lambda}_{ik} \,\, \mbox{ and } \, \, Y_{kj}:=\sum_{q=1}^{l}y_{l-q+1}G_{kj}^{q} \in \tilde{\Lambda}_{kj.}$$
Then
$$X_{ik}Y_{kj}
=\left\{
\begin{array}{ll}
		0 &   if\;m+l-\lambda_{g(k)}\textless 1, \\  \\ \sum\limits_{u=1}^{m+l-\lambda_{g(k)}}\sum\limits_{v=1}^{m+l-\lambda_{g(k)}
        -u+1}x_vy_{m+l-\lambda_{g(k)}-u+1-v+1}G_{ij}^{u} & if\;m+l-\lambda_{g(k)}\ge 1.
\end{array}
\right.$$
\end{Lem}

{\it Proof.} (1) By definition,
\begin{align*}
G_{ik}^pG_{kj}^q & = (\sum_{u=1}^{p}e_{p-u+1,\lambda_{g(k)}-u+1})
(\sum_{v=1}^{q}e_{q-v+1,\lambda_{g(j)}-v+1})
=\sum_{u=1}^{p}\sum_{v=1}^{q}\delta_{\lambda_{g(k)}-u+1,q-v+1}e_{p-u+1,\lambda_{g(j)}-v+1}\\
&=\sum_{v\in V}e_{p+q-\lambda_{g(k)}-v+1,\lambda_{g(j)}-v+1,}
\end{align*}
where $V:=\{v\in [q]\mid 1\le \lambda_{g(k)}-q+v\le p\}$.
If $V=\emptyset$, then $G_{ik}^pG_{kj}^q=0$. Now we take into account the case $V\neq\emptyset$.
Let $v_0\in V$, that is, $1\le v_0\le q$ and $1\le \lambda_{g(k)}-q+v_0\le p$.
It follows from $q\le \theta_{kj}=\min\{\lambda_{g(k)},\lambda_{g(j)}\}\le \lambda_{g(k)}$ that $1\le v_0\le p+q-\lambda_{g(k)}$ and $V=\{v\mid 1\le v\le p+q-\lambda_{g(k)}\}$.
Conversely, if $1\le p+q-\lambda_{g(k)}$, then $1\in V$.
Thus $V\neq \emptyset$ if and only if $ p+q-\lambda_{g(k)}\ge 1$. Therefore, for $V\neq\emptyset$,
$$G_{ik}^pG_{kj}^q=\sum_{v=1}^{p+q-\lambda_{g(k)}}e_{p+q-\lambda_{g(k)}-v+1,\lambda_{g(j)}-v+1}=G_{ij}^{p+q-\lambda_{g(k)}}
$$
This completes the proof of (1).

$(2)$ By definition,
$X_{ik}Y_{kj}=(\sum_{p=1}^{m}x_{m-p+1}G_{ik}^{p})(\sum_{q=1}^{l}y_{l-q+1}G_{kj}^{q})
=\sum_{p=1}^{m}\sum_{q=1}^{l}x_{m-p+1}y_{l-q+1}G_{ik}^{p}G^{q}_{kj.}$
By (1), $G_{ik}^{p}G^{q}_{kj}=0$ for $p+q-\lambda_{g(k)}<1$.
Let $D:=\{(p,q)\mid 1\le p\le m,1\le q\le l, 1\le p+q-\lambda_{g(k)}\}$.
Then $G_{ik}^{p}G^{q}_{kj}=0$ for $(p,q)\not\in D$, and therefore
$$(*) \quad \sum_{p=1}^{m}\sum_{q=1}^{l}x_{m-p+1}y_{l-q+1}G_{ik}^{p}G^{q}_{kj}=\sum_{(p,q)\in D}x_{m-p+1}y_{l-q+1}G_{ik}^{p}G^{q}_{kj.}$$
If $D=\emptyset$, then the summation $(*)$ equals $0$, and therefore $X_{ik}Y_{kj}=0$. If $D\neq\emptyset$, then we pick up an arbitrary element $(p_0,q_0)\in D$, that is, $1\le p_0\le m, 1\le q_0\le l$ and $1\le p_0+q_0-\lambda_{g(k)}$. In this case, $\lambda_{g(k)}-p_0+1\le q_0\le l$ and $\lambda_{g(k)}-l+1\le  p_0\le m$.
Hence $m+l-\lambda_{g(k)}\ge 1 $ and $D=\{(p,q)\mid \lambda_{g(k)}-l+1\le p\le m,\lambda_{g(k)}-p+1\le q\le l\}$.
Conversely, if $ m+l-\lambda_{g(k)}\ge 1$, then $(m,l)\in D$.
Thus $D\neq\emptyset$ if and only if $m+l-\lambda_{g(k)}\ge 1$.
It follows from (1) that, for $m+l-\lambda_{g(k)}\ge 1$,
\begin{align*}
X_{ik}Y_{kj}&=\sum_{p=\lambda_{g(k)}-l+1}^{m}\sum_{q=\lambda_{g(k)}-p+1}^{l}x_{m-p+1}y_{l-q+1}G^{p+q-\lambda_{g(k)}}_{ij}
=\sum_{p=\lambda_{g(k)}-l+1}^{m}\sum_{u=1}^{p+l-\lambda_{g(k)}}x_{m-p+1}y_{p+l-\lambda_{g(k)}-u+1}G^{u}_{ij}\\
&=\sum_{v=1}^{m+l-\lambda_{g(k)}}\sum_{u=1}^{m+l-\lambda_{g(k)}-v+1}x_{v}y_{m+l-\lambda_{g(k)}-u+1-v+1}G^{u}_{ij}
=\sum_{u=1}^{m+l-\lambda_{g(k)}}\sum_{v=1}^{m+l-\lambda_{g(k)}-u+1}x_{v}y_{m+l-\lambda_{g(k)}-u+1-v+1}G^{u}_{ij.}
\end{align*}
Thus (2) follows. $\square$

\medskip
If we write $a\in M_n(R)$ as an $m_s\times m_s$ block matrix $a=(A_{pq})$ with $A_{pq}\in M_{\lambda_{g(i)}\times \lambda_{g(j)}}(R)$, then $a\in f_iS_n(c,R)f_j$ if and only if the $(i,j)$-block $A_{ij}$ satisfies $J_{g(i)}A_{ij}=A_{ij}J_{g(j)}$ and $A_{pq}=0$ for $(p,q) \ne (i,j)$.
Thus there is an  isomorphism of $R$-modules
\begin{align*}
 \tilde{\Lambda}_{ij} &\lra \Lambda_{ij}=f_i S_n(c,R) f_j \\
  a=\sum_{p=1}^{\theta_{ij}} a_p (\sum_{u=1}^{p}e_{p-u+1,\lambda_{g(j)}-u+1}) & \mapsto \sum_{p=1}^{\theta_{ij}} a_p (\sum_{u=1}^{p}e_{n_{g(i)h(i)}-\lambda_{g(i)}+p-u+1,n_{g(j)h(j)}-u+1}),
\end{align*}
induced by the canonical monomorphism of $R$-modules
$$ \varphi_{ij}:M_{\lambda_{g(i)}\times \lambda_{g(j)}}(R) \lra M_n(R),
\;b=\sum_{p=1}^{\lambda_{g(i)}}\sum_{u=1}^{\lambda_{g(j)}}b_{pu}e_{pu}\mapsto \sum_{p=1}^{\lambda_{g(i)}}\sum_{u=1}^{\lambda_{g(j)}}b_{pu}e_{n_{g(i)h(i)}-\lambda_{g(i)}+p,n_{g(j)h(j)}-\lambda_{g(j)}+u},$$
which sends $b=(b_{pq})\in M_{\lambda_{g(i)}\times \lambda_{g(j)}}(R)$ to an $m_s\times m_s$ block matrix in which $b$ is in the $(i,j)$-block and $0$ in all other blocks.
Let
$$F_{ij}^{p}:=(G_{ij}^{p})\varphi_{ij}=	\sum_{u=1}^{p}e_{n_{g(i)h(i)}-\lambda_{g(i)}+p-u+1,n_{g(j)h(j)}-u+1.}$$
Then $\{F_{ij}^{p}\mid 1\le p\le \theta_{ij}\}$ is an $R$-basis of $f_i S_n(c,R) f_j$ and
$\{F_{ij}^{p}\mid 1\le i,j\le m_s,1\le p\le \theta_{ij}\}$ is an $R$-basis of $S_n(c,R)$.

\begin{Lem}\label{multbasis2}
If $1\le i,j,k,l\le m_s,1\le p\le \theta_{ik},1\le q \le \theta_{kj}$, then $$ F_{ik}^pF_{lj}^q=\delta_{kl}(G_{ik}^pG_{kj}^q)\varphi_{ij}=\left\{
	\begin{array}{ll}
		0   & if \; k\neq l \mbox{ or } \;p+q-\lambda_{g(k)}\textless 1,\\
        F_{ij}^{p+q-\lambda_{g(k)}}  & if\; k=l \mbox{ and } p+q-\lambda_{g(k)}\ge 1.
	\end{array} \right.$$
\end{Lem}

{\it Proof.} Clearly, $F_{ik}^p\in f_i\Lambda f_k$ and $F_{lj}^q\in f_l\Lambda f_j$.
If $k\neq l$, then $f_kf_l=0$, and therefore $F_{ik}^pF_{lj}^q=0$.
If $k=l$, then
\begin{align*}
F_{ik}^pF_{kj}^q & = (\sum_{u=1}^{p}e_{n_{g(i)h(i)}-\lambda_{g(i)}+p-u+1,n_{g(k)h(k)}-u+1})
(\sum_{v=1}^{q}e_{n_{g(k)h(k)}-\lambda_{g(k)}+q-v+1,n_{g(j)h(j)}-v+1})\\
&=\sum_{u=1}^{p}\sum_{v=1}^{q}\delta_{n_{g(k)h(k)}-u+1,n_{g(k)h(k)}-\lambda_{g(k)}+q-v+1}e_{n_{g(i)h(i)}-\lambda_{g(i)}+p-u+1,n_{g(j)h(j)}-v+1}\\
&=(\sum_{u=1}^{p}\sum_{v=1}^{q}\delta_{\lambda_{g(k)}-u+1,q-v+1}e_{p-u+1,\lambda_{g(j)}-v+1})\varphi_{ij}\\ &=(G_{ik}^pG_{kj}^q)\varphi_{ij}.
\end{align*}
If $p+q-\lambda_{g(k)}\textless 1$, then $G_{ik}^pG_{kj}^q=0$ by Lemma \ref{multbasis1}(1), and therefore $F_{ik}^pF_{kj}^q=0$.
If $p+q-\lambda_{g(k)}\ge 1$, then it follows from Lemma \ref{multbasis1}(1) that $F_{ik}^pF_{kj}^q=(G_{ik}^pG_{kj}^q)\varphi_{ij}=(G_{ij}^{p+q-\lambda_{g(k)}})\varphi_{ij}=F_{ij}^{p+q-\lambda_{g(k)}}.$
$\square$
	
\begin{Lem}\label{LEM1STEP1} Let $c$ = \emph{diag}$(J_{1}^{b_1},J_{2}^{b_2},\dots,J_{s}^{b_s})\in M_{n}(R)$ be a Jordan-block matrix with Jordan blocks $J_{j}=[r,1,0,\cdots,0]\in M_{\lambda_{j}}(R)$ for $1\leq j\leq s$ where $J_j$ appears $b_j$ times. If $\lambda_1 \textgreater \lambda_2 \textgreater \cdots \textgreater \lambda_s$, then
		
$(1)$ $f_i S_n(c,R) f_i\simeq R[x]/(x^{\lambda_{g(i)}})$ for $1\leq i\leq m_s$.
		
$(2)$ $f_i S_n(c,R) f_j$ is a free $R$-module with an $R$-basis $\{F_{ij}^{p}\mid 1\le p\le \theta_{ij}\}$ of rank $\theta_{ij}$ for $1\le i,j\le m_s$.
		
$(3)$ $S_n(c,R)$ is a free $R$-module with an $R$-basis $\{F_{ij}^{p}\mid 1\le i,j\le m_s,1\le p\le \theta_{ij}\}$ of rank $\sum_{i=1}^{s}(m^2_i-m_{i-1}^2)\lambda_i$.

$(4)$ If $R$ is a local ring and all $b_i=1$, then

{\rm (i)} $\{f_1,f_2,\cdots,f_s\}$ is a complete set of primitive orthogonal idempotents of $S_n(c,R)$.
		
{\rm (ii)} $S_n(c,R)$ is a basic algebra, that is, $S_n(c,R)f_i\ncong S_n(c,R)f_j$ for all $i,j$ with $1\le i\neq j\le s$.
		
{\rm (iii)} $\rad(S_n(c,R))=\{a=(A_{ij})_{s\times s}\in S_n(c,R)\mid A_{ii}=[x_{i1},x_{i2},\dots,x_{i\lambda_i}], \; x_{i1}\in\rad(R), x_{ij}\in R, 1\le i\le s, 2\le j\le \lambda_i\}$. Particularly, if $R$ is a field, then $\rad(S_n(c,R))$ is a free $R$-module with an $R$-basis $\{F_{ij}^{p}\mid 1\le i,j\le s,1\le p\le \theta_{ij}\} \backslash \{F^{\lambda_i}_{ii}\mid 1\le i\le s\}$.

{\rm (iv)} $S_n(c,R)$ is an indecomposable algebra, that is, $1$ and $0$ are the only central idempotents of $S_n(c,R).$
\end{Lem}

{\it Proof.} (1)-(3) are clear. It remains only to prove (4). Since $b_i=1$ for $1\le i\le s$, we have $g(i)=i$ and $\theta_{ij}=\min\{\lambda_i,\lambda_j\}$ for $1\le i,j\le s$.

(i) Since $R$ is a local ring, the only idempotent elements in $R$ are $0$ and $1$. This implies that the only idempotent elements in $R[x]/(x^{\lambda_i})$ are $0$ and $1$, too. It follows from (1) that $f_iS_n(c,R)f_i$ contains only $0$ and $1$ as idempotent elements. Hence $f_i$ is primitive, and therefore all $f_i$ form a complete set of primitive orthogonal idempotents of $S_n(c,R)$.

(ii) By (2), $S_n(c,R)f_i$ is a free $R$-module of rank $i\lambda_i+\sum_{p=i+1}^{s}\lambda_p$ for $1\le i\le s$. For $1\le i< j\le s$, it follows from $\lambda_i\textgreater \lambda_j$ that the $R$-rank of $S_n(c,R)f_i$ is bigger than the $R$-rank of $S_n(c,R)f_j$. Thus the $R$-rank of $S_n(c,R)f_i$ does not equal the $R$-rank o f $S_n(c,R)f_j$ for $1\le i\neq j\le s$. Hence $S_n(c,R)f_i\ncong S_n(c,R)f_j$ for $i\neq j$. This implies that $S_n(c,R)$ is a basic algebra.

(iii) Let $N:=\{a=(A_{ij})_{s\times s}\in S_n(c,R)\mid A_{ii}=[x_{i1}, x_{i2},\dots,x_{i\lambda_i}], \; x_{i1}\in\rad(R), x_{ij}\in R, 1\le i\le s, 2\le j\le \lambda_i\}$. Then $f_iNf_j=f_iS_n(c,R)f_j$ for $1\le i\neq j\le s$.
If $i=j$, then $f_iNf_i=\{([x_{1}, x_{2},\dots,x_{\lambda_i}])\varphi_{ii}\mid x_{1}\in\rad(R), x_{u}\in R, 2\le u\le \lambda_i\}$. By Lemma \ref{lemsemicir}(2), for $x_{1}, x_{2},\dots,x_{\lambda_i}\in R$, we have  $$([x_{1}, x_{2},\dots,x_{\lambda_i}])\varphi_{ii}=(\sum_{p=1}^{\lambda_i}x_{\lambda_i-p+1}G^p_{ii})\varphi_{ii} =\sum_{p=1}^{\lambda_i}x_{\lambda_i-p+1}F_{ii}^p=x_{1}F_{ii}^{\lambda_i}+\sum_{p=1}^{\lambda_i-1}x_{\lambda_i-p+1}F_{ii}^p.$$
If we take $x_1=0$, then we deduce $F^p_{ii}\in f_iNf_i$ for $1\le p<\lambda_i$, and therefore $\{F_{ii}^{p}\mid 1\le p< \lambda_{i}\}\subseteq f_iNf_i$.
If $R$ is a field, then $\rad(R)=0$ and $f_iNf_i$ is an $R$-space with an $R$-basis $\{F_{ii}^{p}\mid 1\le p< \lambda_{i}\}$. In this case, $N$ is an $R$-space with an $R$-basis $\{F_{ij}^{p}\mid 1\le i,j\le s,1\le p\le \theta_{ij}\} \backslash \{F^{\lambda_i}_{ii}\mid 1\le i\le s\}$ by (3).

For $1\leq i\leq s$, $S_n(c,R)f_i$ is indecomposable and projective by (1), and therefore $f_iS_n(c,R)f_i$ is a local ring. It follows from \cite[Proposition 17.19]{AF} that $S_n(c,R)f_i$ is a projective cover of a simple module. This means that the radical of $S_n(c,R)f_i$ is the only maximal submodule of $S_n(c,R)f_i$. If $Nf_i$ is a maximal submodule of $S_n(c,R)f_i$, then $Nf_i$ is the radical of the $S_n(c,R)$-module $S_n(c,R)f_i$. This implies $\rad(S_n(c,R))=N$. So, we need only to show that $Nf_i$ is a maximal submodule of $S_n(c,R)f_i$ for $1\leq i\leq s$.
For this purpose, we first prove that $N f_i$ is an $S_n(c,R)$-submodule of $S_n(c,R)f_i$.
It is an $R$-submodule of $S_n(c,R)f_i$. We have to prove $f_jS_n(c,R)Nf_i\subseteq f_jNf_i\subseteq Nf_i$ for $1\le j\le s$.
For $ j\neq i$, we have $f_jNf_i=f_jS_n(c,R)f_i$. Then $f_jS_n(c,R)Nf_i\subseteq f_jS_n(c,R)f_i=f_jNf_i\subseteq Nf_i$. For $j=i$, we have to show
$f_iS_n(c,R)Nf_i\subseteq Nf_i$.
This is equivalent to saying that $f_iS_n(c,R)f_kNf_i\subseteq f_iNf_i\subseteq Nf_i$ holds for all $1\le k\le s$.
So, suppose
$$a=\sum_{p=1}^{\theta_{ik}}a_{\theta_{ik}-p+1}F_{ik}^p\in f_iS_n(c,R)f_k,\;b=\sum_{q=1}^{\theta_{ik}}b_{\theta_{ik}-q+1}F_{ki}^q\in f_kNf_i,$$
with all $a_i,b_i\in R$, where $\theta_{ik}=\mbox{min}\{\lambda_i,\lambda_k\}$,we show $ab\in f_iNf_i\subseteq Nf_i$. Actually, by Lemma \ref{multbasis1}(2), if $2\theta_{ik}-\lambda_k<1$, then $ab=0 \in f_iNf_i\subseteq Nf_i$. Now, suppose $2\theta_{ik}-\lambda_k\ge 1$.
By definition,
\begin{align*}
ab &=(\sum_{p=1}^{\theta_{ik}}a_{\theta_{ik}-p+1}F_{ik}^p)
(\sum_{q=1}^{\theta_{ik}}b_{\theta_{ik}-q+1}F_{ki}^q)\\
&=\sum_{1\le p,q\le \theta_{ik}}a_{\theta_{ik}-p+1}b_{\theta_{ik}-q+1}(F_{ik}^pF_{ki}^q)
\quad (\mbox{by Lemma }\ref{multbasis2})\\
&=\sum_{1\le p,q\le \theta_{ik}}a_{\theta_{ik}-p+1}b_{\theta_{ik}-q+1}
(G_{ik}^pG_{ki}^q)\varphi_{ii}
\quad (\mbox{ by Lemma \ref{multbasis1}(2) })\\
\end{align*}
\begin{align*}
&= (\sum_{u=1}^{2\theta_{ik}-\lambda_k}
\sum_{v=1}^{2\theta_{ik}-\lambda_k-u+1}a_vb_{2\theta_{ik}-\lambda_k-u+1-v+1}G_{ii}^u)\varphi_{ii.}
\end{align*}
If $i\neq k$, then $2\,\theta_{ik}-\lambda_k=2\,\mbox{min}\{\lambda_i,\lambda_k\}-\lambda_k<\lambda_i$ and $\lambda_i-(2\,\theta_{ik}-\lambda_k)\ge 1$. In this case,  $$ab=([0,\cdots,0,a_1 b_1,a_1 b_2+a_2 b_1,\cdots,\sum_{v=1}^{2\theta_{ik}-\lambda_k}a_vb_{2\theta_{ik}-\lambda_k-v+1}])\varphi_{ii,}$$
where $0$ appears $\lambda_i-(2\,\theta_{ik}-\lambda_k)$ times.
Note that $f_iNf_i=\{([x_{1}, x_{2},\dots,x_{\lambda_i}])\varphi_{ii}\mid x_{1}\in\rad(R), x_{u}\in R, 2\le u\le \lambda_i\}$.
Thus $ab\in f_iNf_i\subseteq Nf_i$.

If $i=k$, then $\theta_{ik}=\mbox{min}\{\lambda_i,\lambda_k\}=\lambda_i$, and therefore  $$ab=(\sum_{u=1}^{\lambda_i}\sum_{v=1}^{\lambda_i-u+1}a_vb_{\lambda_i-u+1-v+1}G_{ii}^u)\varphi_{ii.}$$
that is, $ab$ is an $s\times s$ block matrix with $[a_1 b_1,a_1 b_2+a_2 b_1,\cdots,\sum_{v=1}^{\lambda_i}a_vb_{\lambda_i-v+1}]$ in the $(i,i)$-block of size $\lambda_{i}\times \lambda_{i}$, and $0$ in the $(p,q)$-block of size $\lambda_{p}\times \lambda_{q}$ if $(p,q)\ne (i,i)$.
As $f_iNf_i=\{([x_{1}, x_{2},\dots,x_{\lambda_i}])\varphi_{ii}\mid x_{1}\in\rad(R), x_{u}\in R, 2\le u\le \lambda_i\}$, it follows from $b_1\in \rad(R)$ that $a_1 b_1\in \rad(R)$ and $ab\in f_iNf_i\subseteq Nf_i$.
Hence $Nf_i$ is a submodule of the $S_n(c,R)$-module $S_n(c,R)f_i$.

Now, we show that $Nf_i$ is a maximal submodule of $S_n(c,R)f_i$.
Suppose that $M$ is a submodule of $S_n(c,R)f_i$ with $Nf_i\subsetneq M$.
Since $Nf_i$ is properly contained in $M$, there is an element $y\in M\setminus Nf_i$. Since $y=f_1y+\cdots +f_iy+\cdot +f_sy$ and $f_jS_n(c,R)f_i\subseteq Nf_i$ for $j\ne i$, we deduce $f_iy=f_iyf_i\not\in Nf_i$.
This means that if we write $(f_iyf_i)\varphi_{ii}^{-1}=[y_{1},y_{2},\cdots,y_{\lambda_i}]$, then $y_1\not\in \rad(R)$.
Since $R$ is a local ring, the elements in $R\setminus \rad(R)$ are invertible in $R$.
Thus
\begin{align*}
  F_{ii}^{\lambda_i} & =([1,0,\cdots,0])\varphi_{ii}= \big(y_1^{-1}([y_{1},y_{2},\cdots,y_{\lambda_i}]-[0,y_{2},\cdots,y_{\lambda_i}])\big)\varphi_{ii}\\
  &=y_1^{-1}([y_{1},y_{2},\cdots,y_{\lambda_i}])\varphi_{ii}-y_1^{-1}([0,y_{2},\cdots,y_{\lambda_i}])\varphi_{ii.}
\end{align*}
Thanks to $([y_{1},y_{2},\cdots,y_{\lambda_i}])\varphi_{ii}=f_iyf_i\in M$ and $([0,y_{2},\cdots,y_{\lambda_i}])\varphi_{ii} \in f_iNf_i\subseteq M$, we have $ F_{ii}^{\lambda_i}\in M$.
Moreover, $\{F_{ui}^{p}\mid 1\le u\le s,1\le p\le \theta_{ui}\}\backslash\{F_{ii}^{\lambda_i}\}  \subseteq Nf_i \subset M$. This means that $M$ contains an $R$-basis of $S_n(c,R)f_i$ by (3), and therefore $M=S_n(c,R)f_i$. Hence $Nf_i$ is a maximal submodule of $S_n(c,R)f_i$, and $\rad(S_n(c,R))=N$.

(iv) This follows from the fact that $f_iS_n(c,R)f_j\neq 0$ for all $i,j$ by (3).
$\square$
	
\medskip
Finally, we assume that \textbf{$R$ has no zero-divisors and $c$ is a Jordan-block matrix with different eigenvalues in $Z(R)$.}

In this case we may suppose that $c$ has $t$ distinct eigenvalues in $R$, say $r_1, \cdots, r_t$, and that, for each eigenvalue $r_i$, there are $s_i$ Jordan-blocks $J_{ij}$ of different size $\lambda_{ij}$ with the same eigenvalue $r_i$ for $1\le j\le s_i$, such that $J_{ij}$ appears $b_{ij}$ times in $c$, that is,
$$(\dag\dag)\quad c = \mbox{ diag}(J_{11}^{b_{11}},J_{12}^{b_{12}},\cdots,J_{1s_1}^{b_{1s_1}},J_{21}^{b_{21}}, J_{22}^{b_{22}}\cdots,J_{2s_2}^{b_{2s_2}},\cdots,J_{t1}^{b_{t1}},J_{t2}^{b_{t2}}\cdots,J_{ts_t}^{b_{ts_t}})\in M_n(R),$$
where $J_{ij}=[r_i,1,0,\dots,0]\in M_{\lambda_{ij}}(R)$ and $b_{ij}\ge 1$ for $1\leq j\leq s_i$ and $1\leq i\leq t$. Further, we assume $\lambda_{i1}>\lambda_{i2}>\cdots> \lambda_{is_i}$, $r_i\in Z(R)$, and $r_i\neq r_j$ for $i\neq j$.

Let $\underline{\lambda}_i:=(\lambda_{i1},\lambda_{i2},\cdots,\lambda_{is_i})$ (with a fixed ordering).
The set $\{\underline{\lambda}_i\mid 1\le i\le t\}$ is called the \emph{block type} of $c$. If $t=1$, that is, $c$ is a Jordan-block matrix with the same eigenvalues, then the block type of $c$ just encodes the different sizes of blocks in $c$.
If $d$ is similar to a Jordan-block matrix $c$ by an invertible matrix over $R$, then the block type of $d$ is defined to be the block type of $c$.

We define
$$ n_i:=\sum_{p=1}^{s_i}b_{ip}\lambda_{ip}, \;\tau_0:=0,\;\tau_i:=\sum_{p=1}^{i}n_p,\;\epsilon_i:=\sum_{p=\tau_{i-1}+1}^{\tau_i}e_{pp}\in M_n(R),\;1\leq i\leq t.$$
Note that
$[n]=\{\tau_{q-1}+v\in \mathbb{N}\mid 1\le q\le t,1\le v\le n_q\}$ and that $1=\sum_{i=1}^{t}\epsilon_i$ is a decomposition of $I_n$ into pairwise orthogonal central idempotents in $S_n(c,R)$.
For $1\leq i\leq t$, we define $c_i:=$ diag$(J_{i1}^{b_{i1}},J_{i2}^{b_{i2}},\cdots,J_{is_i}^{b_{is_i}})\in M_{n_i}(R)$. It follows from $r_i\neq r_j$ for $1\leq i\neq j\leq t$ and Lemma \ref{lemsemicir}(2) that $S_n(c,R)$ is isomorphic to $$\mbox{diag}(S_{n_1}(c_1,R),S_{n_2}(c_2,R),\dots,S_{n_t}(c_t,R))$$ as rings. We write these observations as the following lemma for reference.

\begin{Lem}\label{LEMSTEP2} If $R$ has no zero-divisors and $c$ is of the form $(\dag\dag)$, then
		
$(1)$ $1=\sum_{i=1}^{t}\epsilon_i$ is a decomposition of $1$ into pairwise orthogonal central idempotents in $S_n(c,R)$.
		
$(2)$ $S_n(c,R)\simeq\emph{diag}(S_{n_1}(c_1,R),S_{n_2}(c_2,R),\dots,S_{n_t}(c_t,R))$ as rings, that is, $S_n(c,R)$ $\simeq S_{n_1}(c_1,R) \times S_{n_2}(c_2,R) \times \dots\times S_{n_t}(c_t,R).$
\end{Lem}

Thus the study of $S_n(c,R)$ related to a general Jordan-similar matrix $c$ can be reduced to the ones related to Jordan-similar matrices with the same eigenvalues, while the latter cannot be further reduced by Example \ref{ex1}.

For the convenience of the reader, we quote the following elementary fact which will be used frequently in proofs.
\begin{Lem}\label{doublesum}
Suppose that $G$ is an additive group and $n$ is a positive integer. If $a_{pq}\in G$ for $1\le q\le p\le n$, then
$$
\sum_{p=1}^{n}\sum_{q=1}^{p}a_{pq} =\sum_{q=1}^{n}\sum_{p=q}^{n}a_{pq}=\sum_{q=1}^{n}\sum_{u=1}^{n-q+1}a_{u+q-1,q}
=\sum_{u=1}^{n}\sum_{q=1}^{n-u+1}a_{u+q-1,q.}
$$
\end{Lem}

\section{Frobenius extensions\label{FROB}}
This section is devoted to proving Theorem \ref{MAIN1}.

Let $A$ be a unitary ring. If $B$ a subring of $A$ with the same identity, then we say that $B\subseteq A$ is an \emph{extension} of rings. Frobenius extensions, initiated by Kasch, form one of the most prominent instances of extensions of rings. They are a generalization of Frobenius algebras over a field (see \cite{Nakayama1960} and \cite{Kasch1961}) and have played an important role in many aspects of mathematics from representation theory (see \cite{Xi2020},\cite{Xi2019}), knot theory and solutions to Yang-Baxter equitions (see
\cite{Kadison1999}), to topological quantum field theories in lower dimensions and code theory (see \cite{Kock} and \cite{g-wood}). A good introduction to the subject for beginners may be the book  by Kadison (see \cite{Kadison1999}).

\begin{Def}\label{DEFFROB} $(1)$ An extension $B\subseteq A$ of rings is called a \emph{Frobenius extension} if $_BA$ is a finitely generated projective $B$-module and $\Hom_B(_BA,B)\simeq {}_AA_B$ as $A$-$B$-bimodules.
	
	$(2)$ An extension $B\subseteq A$ is said to be \emph{separable} if the multiplication map $A\otimes_BA\ra A, x\otimes y\mapsto xy$, is a split surjective homomorphism of $A$-bimodules, and \emph{split} if the $B$-bimodule ${}_BB_B$ is a direct summand of ${}_BA_B$.
\end{Def}

We need the following properties of Frobenius extensions. For proofs, we refer to \cite[Theorem 1.2, p.3; Corollaries 2.16-17, p.15]{Kadison1996}.

\begin{Lem}\label{FROBLEM}Let $B\subseteq A$ be an extension of rings.
	
	$(1)$ The extension is a Frobenius extension if and only if there exist a $B$-$B$-bimodule homomorphism $E\in \Hom_{B-B}(_BA_B, {}_BB_B)$, and elements $x_i, y_i\in A$, $1\le i\le n,$ such that, for any $a\in A,$
	$$\sum_{i=1}^n x_i\,(y_ia)E = a = \sum_i (ax_i)E\, y_{i.}$$ In this case, $(E, x_i,y_i)$ is called a \emph{Frobenius system} of the extension.
	
	$(2)$ Suppose that $B\subseteq A$ is a Frobenius extension with a Frobenius system $(E,x_i,y_i)$. Then $B\subseteq A$ is split if and only if there exists $d\in C_A(B):=\{a\in A\mid ab=ba \mbox{ for all }\, b\in B\}$ such that $E(d)=1$, and separable if and only if there exists $d\in C_A(B)$ such that $\sum_{i=1}^nx_idy_i=1$.
\end{Lem}

Immediately, we have the basic observations.

\begin{Lem}\label{FROBDP}
If $B_i\subseteq A_i$ is a Frobenius extension of rings for $1\leq i\leq s$, then $B_1\times B_2\times\dots\times B_s\subseteq A_1\times A_2\times \dots \times A_s$ is a Frobenius extension of rings.
\end{Lem}

The following lemma is easy and its proof is left to the reader.

\begin{Lem}\label{FROBISO} Let $f:B\hookrightarrow A$ and $g:C \hookrightarrow A$ be extensions of rings. Assume that there exist ring isomorphisms $\varphi:B\ra C$ and $\psi:A\ra A$ such that $f\psi=\varphi g$. Then $B\subseteq A$ is a Frobenius extension if and only if $C\subseteq A$ is a Frobenius extension. Moreover, $B\subseteq A$ is separable (or split) if and only if $C\subseteq A$ is separable (or split).
\end{Lem}

Consequently, we have the following lemma from Lemmas \ref{FROBISO} and \ref{PROPCMA}(1).

\begin{Lem}\label{FROBISOA} Suppose that two matrices $c$ and $d$ in $M_n(R)$ are similar. Then $S_n(c,R)\subseteq M_n(R)$ is a Frobenius (or separable, or split) extension if and only if so is $S_n(d,R)\subseteq M_n(R)$.
\end{Lem}

Our consideration on centralizer matrix algebras is divided into two cases. First, we consider the centralizers of invertible matrices.

Let $G$ be a subgroup of $GL_n(R)$. An element $i\in [n]$ is called a $G$-\emph{free point} if $g_{ii}=0$ for all $g=(g_{ij})\in G\setminus \{I_n\}$. If $G=\langle \sigma \rangle$ is a cyclic group generated by $\sigma\in GL_n(R)$, then a $G$-free point will simply be called a $\sigma$-\emph{free point}.

\begin{Theo} \label{FROBTHM}
If $R$ is a unitary ring and $G$ is a finite subgroup of $GL_n(R)$ with a $G$-free point, then $S_{n}(G,R)\subseteq M_n(R)$ is a separable Frobenius extension.
\end{Theo}
{\it Proof.} For any $a\in M_n(R)$, we can check $\sum_{g\in G}gag^{-1}\in S_n(G,R)$ and define
a map $$E:M_n(R)\lra S_n(G,R),\; a\mapsto \sum_{g\in G}gag^{-1}.$$

Let $i\in [n]$ be a $G$-free point, and let $x_j:=e_{ji}$ and $y:=e_{ij}\in M_n(R)$ for $j\in [n]$. Then
we prove that $(E,x_{j},y_{j})$ is a Frobenius system, that is, the following assertions (i) and (ii) hold.

(i) $E$ is a homomorphism of $S_n(G,R)$-$S_n(G,R)$-bimodules. In fact, $(a+b)E=(a)E+(b)E$. Moreover, for $c\in S_n(G,R)$ and $a\in M_n(R)$, it follows from $cg=gc$ for all $g\in G$ that  $$(ca)E=\sum_{g\in G}g(ca)g^{-1}=c(\sum_{g\in G}gag^{-1})=c\; (a)E.$$
Similarly, $(ac)E=(a)E\, c$.

(ii) For any $a\in M_n(R)$, $\sum_{j=1}^nx_j(y_ja)E=a$ and $\sum_{j=1}^n(ax_j)Ey_j=a$. Indeed, we have $\sum_{j=1}^nx_jy_j=\sum_{j}e_{ji}e_{ij}=\sum_{j=1}^ne_{jj}=I_n$ and, for $g=(g_{pq})\in G\setminus \{I_n\}$, $x_jgy_j=e_{ji}ge_{ij}=e_{jj}(g_{ii})=0$ since $i$ is a $G$-free point. Thus $$\sum_{j=1}^nx_j(y_ja)E=\sum_{j=1}^nx_j\sum_{g\in G}g(y_ja)g^{-1}=\sum_{j=1}^n(x_jy_j)a+\sum_{j=1}^n\sum_{g\neq I_n}(x_jgy_j)ah^{-1}=\sum_{j=1}^ne_{jj}a+0=a.$$
Similarly, $\sum_{j=1}^n(ax_j)Ey_j=a$.

Thus $(E,x_j,y_j)$ is a Frobenius system and $S_n(G,R)\subseteq M_n(R)$ is a Frobenius extension by Lemma \ref{FROBLEM}(1). Further, since $\sum_{j=1}^nx_jy_j=I_n$, the Frobenius extension is separable by Lemma \ref{FROBLEM}(2).
$\square$

As a consequence of Theorem \ref{FROBTHM}, we have the corollary.

\begin{Koro} Let $G$ be a finite subgroup of $GL_n(R)$ with a $G$-free point. If $|G|$ is invertible in $R$, then

$(1)$ $S_{n}(G,R)\subseteq M_n(R)$ is a split Frobenius extension, and global and dominant dimensions of $S_n(G,R)$ are the same as the ones of $R$, respectively.

$(2)$ $S_n(G,R)$ is semisimple if $R$ is semisimple.
\end{Koro}

{\it Proof.} (1) Since $$(\frac{1}{|G|})E=\sum_{h\in G}h\frac{1}{|G|}h^{-1}
=\frac{1}{|G|}\sum_{h\in G}hh^{-1}=1\in S_n(G,R),$$
the extension is split by Lemma \ref{FROBLEM}(2). The statement on global dimensions follows from the fact that the extension is a split, separable Frobenius extension. In the case, the global dimension of $S_n(G,R)$ equals the one of $M_n(R)$ (see\cite[p.14]{Kadison1996}), and the latter equals the global dimension of $R$ since $R$ and $M_n(R)$ are Morita equivalent. For dominant dimensions, the statement follows from \cite{Xi2019} or \cite[p.91]{Tachikawa}, where the definition of dominant dimensions can also be founded.

(2) If $R$ is a semisimple ring and $|G|$ is invertible in $R$, then $S_n(G,R)$ is semisimple. This follows from \cite[Theorem 1.14]{Montgomery1980} which says that, for a finite group acting on a ring $R$, the Jacobson radical of the fixed ring is the intersection of the Jacobson radical of $R$ with the fixed ring itself if the order of $G$ is invertible in $R$. $\square$

\medskip
Now, we apply Theorem \ref{MAIN1} to the centralizers of permutation matrices. To state our result precisely, we first introduce a few notions.

For a natural number $n$, we denote by $\Sigma_{n}$ the symmetric group of all permutations on $[n]$. Any $\sigma\in\Sigma_n$ can be written as a product of disjoint cycles, say $\sigma=\sigma_1\sigma_2\dots\sigma_s$, where $\sigma_i$ is a $\lambda_i$-cycle. Here, $\lambda_i$ is allowed to be $1$, and $\lambda=(\lambda_1,\lambda_2,\dots,\lambda_s)$ is called the cycle type of $\sigma$. Clearly, the order of $\sigma$ is the least common multiple of $\lambda_i's$, denoted by lcm$(\lambda_1,\lambda_2,\dots,\lambda_s)$. Let $c_n(\sigma):=e_{1,(1)\sigma}+e_{2,(2)\sigma}+\cdots+e_{n,(n)\sigma}$ be the permutation matrix in $M_n(R)$ corresponding to $\sigma$. If the subscript $n$ in $c_n(\sigma)$ is clear from the context, we simply write $c_{\sigma}$ for $c_n(\sigma)$. We have $c_{\sigma}c_{\tau}=c_{\sigma\tau}$, $c_{\sigma}'=c_{\sigma^{-1}}=(c_{\sigma})^{-1}$ and $(c_{\sigma})_{ij}=\delta_{(i)\sigma,j}$ for $i,j\in [n]$. Recall that $a'$ denotes the transpose of the matrix $a$ and $\delta_{ij}$ is the Kronecker symbol.

Let $G$ be a subgroup of $\Sigma_n$, and let $P_G:=\{c_{\sigma}\mid\sigma\in G\}$ is the set of all permutation matrices induced by the elements of $G$. Then $P_G$ is clearly a finite subgroup of $GL_n(R)$ isomorphic to $G$.
For $\sigma\in G$ and $i\in [n]$, $(c_{\sigma})_{ii}=\delta_{(i)\sigma,i}$. Thus $(c_{\sigma})_{ii}=1_{R}$ if and only if $(i)\sigma=i$. Hence $i$ is a $P_G$-free point if and only if $\{\sigma\in G\mid(i)\sigma=i\}={1_G}$, that is, $i$ is a $P_G$-free point if and only if the stabilizer $st_{G}(i)$ of $i$ under $G$ is trivial. In this case, we also say that $i$ is a $G$-free point. If $G=\langle\sigma\rangle$, then $G$-free points will simply be called $\sigma$-free points. If $X_i$ is the content of $\sigma_i$ for $1\le i\le s$, then $X_i$ forms a $G$-orbit, $|X_i|=\lambda_i$ and $[n]=\bigcup^{s}_{i=1}X_i$. Moreover, $j\in X_i$ is a $\sigma$-free point if and only if the order of $\sigma$ is just $\lambda_i$. This implies that there is a $\sigma$-free point in $[n]$ if and only if there is some $\lambda_j$ such that $\lambda_i|\lambda_j$ for all $1\leq i\leq s$.

 Clearly, $G$ also acts on $M_n(R)$ by
 $$M_n(R)\times G\lra M_n(R),\; (a_{ij})^{\sigma}:=(a_{(i)\sigma,(j)\sigma}),\; (a_{ij})\in M_n(R),\; \sigma\in G.$$
Since $c_{\sigma}ac_{\sigma}^{-1}=(\sum_{i=1}^{n}e_{i,(i)\sigma})(\sum_{p,q}a_{pq}e_{pq}) (\sum_{j=1}^{n}e_{j,(j)\sigma^{-1}})=\sum_{ij}a_{(i)\sigma,(j)\sigma}e_{ij}=a^{\sigma}$ for $a=(a_{pq})\in M_n(R)$, we have $S_n(P_G,R)=\{a\in M_n(R)\mid a^{\sigma}=a, \sigma\in G\}$ = $M_n(R)^G$, the fixed ring of $G$ in $M_n(R)$. For brevity, we write $S_n(G,R)$ for
$S_n(P_G,R)$, and $S_n(\sigma,R)$ for $S_n(\langle\sigma\rangle,R)$ if $\sigma\in \Sigma_n$.

\begin{Koro}\label{PERMFROB} Let $R$ be a ring and $G$ be a subgroup of $\Sigma_n$ with a $G$-free point. Then

$(1)$ $S_n(G,R)\subseteq M_n(R)$ is a separable Frobenius extension.

$(2)$ If $|G|$ is invertible in a ring $R$, then

{\rm (i)} $S_{n}(G,R)\subseteq M_n(R)$ is a split Frobenius extension, and global and dominant dimensions of $S_n(G,R)$ are the same as the ones of $R$, respectively.

{\rm (ii)} $S_n(G,R)$ is semisimple if $R$ is semisimple.\end{Koro}

Note that the condition on $G=\langle \sigma\rangle$ in Corollary \ref{PERMFROB} is satisfied for $\sigma:i\mapsto n+1-i$. Thus Corollary \ref{PERMFROB} extends the first statement of \cite[Theorem (1), p.318]{Xi2020}. 

\medskip
Having considered the centralizers of invertible matrices in Theoem \ref{MAIN1}, we next investigate the centralizers of not necessarily invertible matrices.

Recall that a matrix in $M_n(R)$ is called a Jordan-similar matrix if it is similar to a Jordan-block matrix by a matrix in $GL_n(R)$. As is known, every square matrix over an algebraically closed field is a Jordan-similar matrix.

Now, let $c$ be a Jordan-block matrix as in $(\dag)$, say
$$c= {\rm diag}(J_{1}^{b_1},J_{2}^{b_2},\cdots,J_{s}^{b_s})\in M_{n}(R),$$
with Jordan blocks $J_{j}=[r,1,0,\cdots,0]\in M_{\lambda_{j}}(R)$ for $1\leq j\leq s$ and $r\in Z(R)$, where $J_j$ appears $b_j$ times and $\lambda_1 \textgreater \lambda_2 \textgreater \cdots \textgreater \lambda_s$.

For $1\le i, j\le m_s$, let $\rho_{ij}:=\lambda_{g(i)}+\lambda_{g(j)}-\lambda_1$. Further, we define a map $E_{ij}: M_{\lambda_{g(i)}\times \lambda_{g(j)}}(R)\ra\tilde{\Lambda}_{ij}$ of $R$-modules as follows. If $\rho_{ij}\le 0$, we set $E_{ij}=0$. If  $\rho_{ij}>0$, we define
$$
E_{ij}: M_{\lambda_{g(i)}\times \lambda_{g(j)}}(R)\lra\tilde{\Lambda}_{ij}, \;
a:=(a_{kl})\mapsto
\sum_{p=1}^{\rho_{ij}}\sum_{u=1}^{\rho_{ij}-p+1}a_{\lambda_{g(i)}-u+1,\rho_{ij}-p+1-u+1}G_{ij}^{p}\in \tilde{\Lambda}_{ij.}
$$
Then $E_{ij}$ is a homomorphism of $R$-modules.
Next, we extend $E_{ij}$ to a map $\bar{E}$ from $M_n(R)$ to $\Lambda$ by setting
$$\bar{E}:M_n(R)\lra \Lambda,\; (A_{ij})_{m_s\times m_s}\mapsto \big((A_{ij})E_{ij}\big)=\sum_{1\le i,j\le m_s}\big((A_{ij})E_{ij}\big)\varphi_{ij}, $$
where $A_{ij}\in M_{\lambda_{g(i)}\times \lambda_{g(j)}}(R)$ for $1\leq i,j\leq m_s$.
This map has the property.

\begin{Lem} \label{stepone}
Let
$c= {\rm diag}(J_{1}^{b_1},J_{2}^{b_2},\cdots,J_{s}^{b_s})\in M_{n}(R)$ be a Jordan-block matrix with
Jordan blocks $J_{j}=[r,1,0,\cdots,0]\in M_{\lambda_{j}}(R)$ for $1\leq j\leq s$ and $r\in Z(R)$, where $J_j$ appears $b_j$ times and $\lambda_1 \textgreater \lambda_2 \textgreater \cdots \textgreater \lambda_s$. Then

$(1)$ $\bar{E}$ is a homomorphism of $\Lambda$-$\Lambda$-bimodules.

$(2)$ If $a=(a_{uv})\in M_n(R)$ and $1\le i\le n$, then

$\quad \rm{(i)}$ $e_{11}(e_{\lambda_1,i}a)\bar{E}=e_{1i}a$, that is, $e_{11}(\sum_{p=1}^{n}a_{ip}e_{\lambda_1,p})\bar{E}=\sum_{p=1}^{n}a_{ip}e_{1p}.$

$\quad \rm{(ii)}$ $(ae_{i1})\bar{E}e_{\lambda_1,1}=ae_{i1}$, that is, $(\sum_{p=1}^{n}a_{pi}e_{p1})\bar{E}e_{\lambda_1,1}=\sum_{p=1}^{n}a_{pi}e_{p1}$.

$(3)$ If there exists an $n\times n$ matrix $z\in M_n(R)$ such that $(z)\bar{E}=1$ and $za=az$ for all $a\in \Lambda$, then $c=rI_n$.
\end{Lem}

{\it Proof.} (1) The map $\bar{E}$ is additive because each $E_{ij}$ is additive. Further, we prove $$(xa)\bar{E}=x(a)\bar{E} \mbox{  and } (ax)\bar{E}=(a)\bar{E}x \mbox{ for any } x\in\Lambda \mbox{ and } a\in M_n(R).$$ This is equivalent to saying that $(X_{ik}A_{kj})E_{ij}=X_{ik}(A_{kj})E_{kj}$ for $1\leq i,k,j\leq m_s,\;X_{ik}\in \tilde{\Lambda}_{ik}$ and $ A_{kj}\in M_{\lambda_{g(k)}\times \lambda_{g(j)}}(R)$.

Indeed, let $l:=\mbox{max}\{i,k\}$. Then $\lambda_{g(l)}=\min\{\lambda_{g(i)},\lambda_{g(k)}\}=\theta_{ik}$ and $\rho_{lj}\le \rho_{ij}$.
We may write $X_{ik}$ = $\sum_{p=1}^{\lambda_{g(l)}}x_{\lambda_{g(l)}-p+1}G_{ij}^{p}\in \tilde{\Lambda}_{ik}$ and $ A_{kj}=(a_{uv})\in M_{\lambda_{g(k)}\times \lambda_{g(j)}}(R)$.
Then
\begin{align*}	X_{ik}A_{kj}&=(\sum_{p=1}^{\lambda_{g(l)}}x_{\lambda_{g(l)}-p+1}\sum_{w=1}^{p}e_{p-w+1,\lambda_{g(k)}-w+1})(
\sum_{q=1}^{\lambda_{g(k)}}\sum_{v=1}^{\lambda_{g(j)}}a_{qv}e_{qv})  =\sum_{p=1}^{\lambda_{g(l)}}\sum_{w=1}^{p}\sum_{q=1}^{\lambda_{g(k)}}\sum_{v=1}^{\lambda_{g(j)}}x_{\lambda_{g(l)}-p+1}a_{qv}\delta_{\lambda_{g(k)}-w+1,q}e_{p-w+1,v}\\
&=\sum_{v=1}^{\lambda_{g(j)}}\sum_{p=1}^{\lambda_{g(l)}}\sum_{w=1}^{p}x_{\lambda_{g(l)}-p+1}a_{\lambda_{g(k)}-w+1,v}e_{p-w+1,v} \quad (\mbox{ by Lemma }\ref{doublesum})\\
&= \sum_{v=1}^{\lambda_{g(j)}}\sum_{u=1}^{\lambda_{g(l)}}\sum_{w=1}^{\lambda_{g(l)}-u+1}x_{\lambda_{g(l)}-(u+w-1)+1}a_{\lambda_{g(k)}-w+1,v}e_{uv,}
\end{align*}
We write $X_{ik}A_{kj}$ as a $\lambda_{g(i)}\times \lambda_{g(j)}$ matrix $(b_{uv})$ with the $(u,v)$-entry $b_{uv} = \sum_{w=1}^{\lambda_{g(l)}-u+1}x_{\lambda_{g(l)}-u+1-w+1}a_{\lambda_{g(k)}-w+1,v}$ for $1\le u\le \lambda_{g(l)}$ and $1\le v\leq \lambda_{g(j)}$, and other entries $0$.

Next, we consider $(X_{ik}A_{kj})E_{ij}$. By definition, if $\rho_{ij}\le 0$, then $(X_{ik}A_{kj})E_{ij}=0$.
Suppose $\rho_{ij}>0$.
Then
$$(X_{ik}A_{kj})E_{ij} = \sum_{p=1}^{\rho_{ij}}\sum_{t=1}^{\rho_{ij}-p+1}b_{\lambda_{g(i)}-t+1,\rho_{ij}-p+1-t+1}G_{ij}^{p}.$$
Note that $b_{uv}=0$ for $u>\lambda_{g(l)}$. Let $D:=\{(p,t)\mid 1\le p\le \rho_{ij},1\le t\le \rho_{ij}-p+1,\lambda_{g(i)}-t+1\le \lambda_{g(l)}\}$. Then $b_{\lambda_{g(i)}-t+1,\rho_{ij}-p+1-t+1}=0$ for $(p,t)\not\in D$, and therefore
$$(X_{ik}A_{kj})E_{ij} = \sum_{p=1}^{\rho_{ij}}\sum_{t=1}^{\rho_{ij}-p+1}b_{\lambda_{g(i)}-t+1,\rho_{ij}-p+1-t+1}G_{ij}^{p} = \sum_{(p,t)\in D}b_{\lambda_{g(i)}-t+1,\rho_{ij}-p+1-t+1}G^{p}_{ij,}
$$
If $D=\emptyset$, then $(X_{ik}A_{kj})E_{ij}=0$. Now we take into account the case $D\neq\emptyset$.
Let $(p_0,t_0)\in D$, that is, $1\le p_0\le \rho_{ij},1\le t_0\le \rho_{ij}-p_0+1,\lambda_{g(i)}-t_0+1\le \lambda_{g(l)}$.
It follows from $\lambda_{g(i)}\ge \lambda_{g(l)}$ that $\lambda_{g(i)}-\lambda_{g(l)}+1\le t_0\le \rho_{ij}-p_0+1$.
Then $p_0\le \rho_{ij}-\lambda_{g(i)}+\lambda_{g(l)}=\rho_{lj}$.
Hence $1\le p_0 \le \rho_{lj}$ and $D=\{(p,t)\mid 1\le p\le \rho_{lj},\lambda_{g(i)}-\lambda_{g(l)}+1\le t\le \rho_{ij}-p+1\}$. Conversely, if $\rho_{lj}\ge 1$, then $(1,\rho_{ij})\in D$. Thus $D\neq\emptyset$ if and only if $\rho_{lj}\ge 1.$ So, if $D\neq\emptyset$, then
\begin{align*}
(X_{ik}A_{kj})E_{ij}&= \sum_{(p,t)\in D}b_{\lambda_{g(i)}-t+1,\rho_{ij}-p+1-t+1}G^{p}_{ij}
=\sum_{p=1}^{\rho_{lj}}\sum_{t=\lambda_{g(i)}-\lambda_{g(l)}+1}^{\rho_{ij}-p+1}b_{\lambda_{g(i)}-t+1,\rho_{ij}-p+1-t+1}G_{ij}^{p}\\
&=\sum_{p=1}^{\rho_{lj}}\sum_{q=1}^{\rho_{lj}-p+1}b_{\lambda_{g(l)}-q+1,\rho_{lj}-p+1-q+1}G^{p}_{ij}\\
&=\sum_{p=1}^{\rho_{lj}}\sum_{q=1}^{\rho_{lj}-p+1}\sum_{w=1}^{q}x_{q-w+1}a_{\lambda_{g(k)}-w+1,\rho_{lj}-p+1-q+1}G^{p}_{ij} \quad (\mbox{ by Lemma }\ref{doublesum})\\
&=\sum_{p=1}^{\rho_{lj}}\sum_{v=1}^{\rho_{lj}-p+1}\sum_{w=1}^{\rho_{lj}-p+1-v+1}x_{v}a_{\lambda_{g(k)}-w+1,\rho_{lj}-p+1-v+1-w+1}G^{p}_{ij.}
\end{align*}
Note that we always have $\rho_{lj} \le \rho_{ij}$, and therefore the following holds always
$$(X_{ik}A_{kj})E_{ij}=\left\{
	\begin{array}{ll}
0  & if\;\rho_{lj} < 1,\\ \sum\limits_{p=1}^{\rho_{lj}}\sum\limits_{v=1}^{\rho_{lj}-p+1}
\sum\limits_{w=1}^{\rho_{lj}-p+1-v+1}x_{v}a_{\lambda_{g(k)}-w+1,\rho_{lj}-p+1-v+1-w+1}G^{p}_{ij} & if\; \rho_{lj} \ge 1.
	\end{array} \right.$$

It remains to consider $X_{ik}\,(A_{kj})E_{kj}$. Firstly, $(A_{kj})E_{kj}=0$ for $\rho_{kj}\le 0$ by definition.
Assume $\rho_{kj}>0$. Then
$$(A_{kj})E_{kj}=\sum_{p=1}^{\rho_{kj}}(\sum_{u=1}^{\rho_{kj}-p+1}a_{\lambda_{g(k)}-u+1,
\rho_{kj}-p+1-u+1})G_{kj.}^{p}$$
Note that $X_{ik}=\sum_{p=1}^{\lambda_{g(l)}}x_{\lambda_{g(l)}-p+1}G_{ik}^{p} \in \tilde{\Lambda}_{ik},(A_{kj})E_{kj}\in \tilde{\Lambda}_{kj}$ and $\lambda_{g(l)}+\rho_{kj}-\lambda_{g(k)}=\lambda_{g(l)}+\lambda_{g(k)}+
\lambda_{g(j)}-\lambda_1-\lambda_{g(k)}=\rho_{lj}.$ According to Lemma \ref{multbasis1}(2), if $\rho_{lj}<1$, then $X_{ik}(A_{kj})E_{kj}=0$; if $\rho_{lj}\ge 1$, then
$$X_{ik}(A_{kj})E_{kj}=\sum_{p=1}^{\rho_{lj}}\sum_{v=1}^{\rho_{lj}-p+1}\sum_{w=1}^{(\rho_{lj}-p+1)-v+1}x_v a_{\lambda_{g(k)}-w+1,
(\rho_{lj}-p+1)-v+1-w+1}G^p_{ij.}$$
Since we always have $\rho_{lj}\le \rho_{kj}$, it follows that
$$X_{ik}(A_{kj})E_{kj}=\left\{
	\begin{array}{ll}
0    & if\;\rho_{lj} < 1, \\		\sum\limits_{p=1}^{\rho_{lj}}\sum\limits_{v=1}^{\rho_{lj}-p+1}\sum\limits_{w=1}^{(\rho_{lj}-p+1)-v+1}
x_v a_{\lambda_{g(k)}-w+1,
(\rho_{lj}-p+1)-v+1-w+1}G^p_{ij}  & if\; \rho_{lj} \ge 1.
	\end{array} \right.$$
Hence $(X_{ik}A_{kj})E_{ij}=X_{ik}(A_{kj})E_{kj}$ holds for all $1\le i,k,j\le m_s,\;X_{ik}\in \tilde{\Lambda}_{ik}$ and $ A_{kj}\in M_{\lambda_{g(k)}\times \lambda_{g(j)}}(R)$.
Thus $(xa)\bar{E}=x(a)\bar{E}$  for any $x\in\Lambda$ and $a\in M_n(R).$
This shows that $\bar{E}$ is a homomorphism of $\Lambda$-modules. Similarly, we show that $\bar{E}$ is a homomorphism of right $\Lambda$-modules.
Thus $\bar{E}$ is a homomorphism of $\Lambda$-$\Lambda$-bimodules.

(2) Now, we prove that the equality $e_{11}(e_{\lambda_1,i}a)\bar{E}=e_{1i}a$ holds for $1\le i\le n$.
In fact, the matrix $e_{\lambda_1,i}a$ has the $\lambda_1$-th row equal to the $i$-th row of $a$, and the other rows equal to $0$, that is, $$e_{\lambda_1,i}a=\sum_{p=1}^{n}a_{ip}e_{\lambda_1,p}=\sum_{j=1}^{m_s}(\sum_{v=1}^{\lambda_{g(j)}}a_{i,n_{g(j)h(j)}-\lambda_{g(j)}+v}e_{\lambda_1v})\varphi_{1j,}$$
So $e_{\lambda_1,i}a$ can be written as an $m_s\times m_s$ block matrix with $\sum_{v=1}^{\lambda_{g(j)}}a_{i,n_{g(j)h(j)}-\lambda_{g(j)}+v}e_{\lambda_1v}$ in the $(1,j)$-block for $1\le j\le m_s$ and $0$ in the other blocks.
Then $(e_{\lambda_1,i}a)\bar{E}=\sum_{j=1}^{m_s}((\sum_{v=1}^{\lambda_{g(j)}}a_{i,n_{g(j)h(j)}-\lambda_{g(j)}+v}e_{\lambda_1v})E_{1j})\varphi_{1j}$
and
\begin{align*}
e_{11}(e_{\lambda_1,i}a)\bar{E}
&=\sum_{j=1}^{m_s}(e_{11}(\sum_{v=1}^{\lambda_{g(j)}}a_{i,n_{g(j)h(j)}-\lambda_{g(j)}+v}e_{\lambda_1v})E_{1j})\varphi_{1j}\quad (\mbox{ by the definition of } E_{1j})\\
&=\sum_{j=1}^{m_s}
(\sum_{p=1}^{\lambda_{g(j)}}a_{i,n_{g(j)h(j)}-p+1}e_{1,\lambda_{g(j)}-p+1})\varphi_{1j}\\
&=\sum_{j=1}^{m_s}(\sum_{w=1}^{\lambda_{g(j)}}a_{i,n_{g(j)h(j)}-\lambda_{g(j)}+w}e_{1w})\varphi_{1j}
\quad(\mbox{ by the definition of } \varphi_{1j})\\
&=\sum_{j=1}^{m_s}\sum_{w=1}^{\lambda_{g(j)}}a_{i,n_{g(j)h(j)}-\lambda_{g(j)}+w}e_{1,n_{g(j)h(j)}-\lambda_{g(j)}+w}\\
&=\sum_{p=1}^{n}a_{ip}e_{1p}\\
&=e_{1i}a.
\end{align*}
The second last equality is due to $[n]=\{n_{g(j)h(j)}-\lambda_{g(j)}+w \mid  1\le j\le m_s,1\le w\le \lambda_{g(j)} \}$.
Similarly, we show $(ae_{i1})\bar{E}e_{\lambda_1,1}=ae_{i1}$ for $1\le i\le n$.

(3) We write $z\in M_n(R)$ as an $m_s\times m_s$ block matrix $z=(Z_{ij})$ with $Z_{ij}\in M_{\lambda_{g(i)}\times \lambda_{g(j)}}$ for $1\le i,j\le m_s$.
For $1\le i\le m_s$, we write $Z_{ii}=(z_{uv})\in M_{\lambda_{g(i)}}(R)$.
It follows from $(Z)E=1$ that
$(Z_{ii})E_{ii}=I_{\lambda_{g(i)}}$, that is,
$$ (*)\quad
(Z_{ii})E_{ii}=\sum_{p=1}^{\rho_{ii}}(\sum_{u=1}^{\rho_{ii}-p+1}z_{\lambda_{g(i)}-u+1,
\rho_{ii}-p+1-u+1})G_{ii}^{p}
=I_{\lambda_{g(i)}}= G_{ii}^{\lambda_{g(i)}}.$$
Since $\rho_{ii}=2\lambda_{g(i)}-\lambda_1\le \lambda_{g(i)}$ and $\{G_{ii}^{p}\mid 1\le p\le \rho_{ii}\}$ is an $R$-basis of $\tilde{\Lambda}_{ii}$, we obtain $\rho_{ii}=\lambda_{g(i)}$, and therefore $\lambda_{g(i)}=\lambda_1$ for $1\le i \le m_s$. Particularly, it follows from $g(m_s)=s$ that $\lambda_s=\lambda_{g(m_s)}=\lambda_1$. This implies $s=1$ by our assumption on $\lambda_i$ for $1\le i\le s$. Hence $c$ is a block matrix with all blocks of the same size. Moreover, by comparing the coefficients of $G_{ii}^{\lambda_1}$ in $(*)$, we deduce
$$z_{\lambda_1,1}=1_R.$$

Next, we show that each of these blocks is in fact a $1\times 1$ matrix. For an $m_1\times m_1$ block matrix $a=(A_{ij})\in \Lambda$ with $A_{ij}\in \tilde{\Lambda}_{ij}$, $za=(\sum_{p=1}^{m_1}Z_{ip}A_{pj})$ and $az=(\sum_{p=1}^{m_1}A_{ip}Z_{pj})$.
Thus the condition $za=az$ for all $a\in \Lambda$ is equivalent to saying that $\sum_{p=1}^{m_1}Z_{ip}A_{pj}=\sum_{p=1}^{m_1}A_{ip}Z_{pj}$ holds for all $A_{pj}\in \tilde{\Lambda}_{pj},A_{jp}\in \tilde{\Lambda}_{jp}, 1\le p,i,j\le m_1$.
In particular, for $1\le i\le m_1$, if
$A_{ip}=0$ and $A_{pi}=0$ for $1\le p\neq i\le m_1$,
we get $Z_{ii}A_{ii}=A_{ii}Z_{ii}$ for all $A_{ii}\in \tilde{\Lambda}_{ii}$.
Now, we write $A_{ii}=[a_1,a_2,\cdots,a_{\lambda_1}]=\sum_{p=1}^{\lambda_1}\sum_{q=1}^{p}a_{\lambda_1-p+1}e_{p-q+1,\lambda_1-q+1}
.$
Then
$$A_{ii}Z_{ii}=\sum_{1\le u,v\le \lambda_1}\sum_{k=u}^{\lambda_1}a_{k-u+1}z_{kv}e_{uv,}\quad \mbox{ and }\quad Z_{ii}A_{ii}=\sum_{1\le u,w\le \lambda_1}\sum_{p=\lambda_1-w+1}^{\lambda_1}z_{u,p-(\lambda_1-w)}a_{\lambda_1-p+1}e_{uw.}$$
Suppose $\lambda_1\ge 2$. Then $(A_{ii}Z_{ii})_{\lambda_1-1, 1}=a_2z_{\lambda_1 1}+a_1z_{\lambda_1-1,1}$ and $(Z_{ii}A_{ii})_{\lambda_1-1, 1}=z_{\lambda_1-1,1}a_1$.
Specially, if $a_1=a_2=1$, then it follows from $Z_{ii}A_{ii}=A_{ii}Z_{ii}$ that $z_{\lambda_1 1}=0$. This contradicts to $z_{\lambda_1 1}=1$.
Thus $\lambda_1=1$ and each block is a $1\times 1$ matrix. This means $c=rI_n$.
$\square$

\begin{Theo}\label{STEPTWO} Suppose that $R$ has no zero-divisors. Let $c\in M_n(R)$ be a Jordan-similar matrix with all eigenvalues in $Z(R)$. Then

$(1)$  $S_n(c,R)\subseteq M_n(R)$ is a separable Frobenius extension.

$(2)$ $S_n(c,R)\subseteq M_n(R)$ is a split extension if and only if $c$ is similar to ${\rm diag}(r_1I_{n_1},r_2I_{n_2},\cdots,r_tI_{n_t})$ with $r_i\neq r_j$ for $1\le i\neq j\le t$, $n_j\ge1$ and $\sum_{j=1}^tn_j=n.$

$(3)$ $S_n(c,R)$ is semisimple if and only if $R$ is semisimple and $c$ is similar to ${\rm diag}(r_1I_{n_1},r_2I_{n_2},\cdots,r_tI_{n_t})$, where $r_i\neq r_j$ for $1\le i\neq j\le t$, $n_j\ge1$ and $\sum_{j=1}^tn_j=n.$
\end{Theo}

{\it Proof.} By Lemma \ref{FROBISOA}, we may assume $c=$ diag$(c_1,c_2,\dots,c_t)$ as in $(\dag\dag)$. Recall that $\tau_0:=0,\; \tau_i:=\sum_{p=1}^{i}n_p$ and $n_i:=\sum_{p=1}^{s_i}b_{ip}\lambda_{ip}$ for $1\le i\le t$ with $\sum_{p=1}^{t}\sum_{v=1}^{n_p}\tau_{p-1}+v=n$ (see the end of Section \ref{CMA}).
The ring $S_n(c,R)$ is isomorphic to the ring diag$(S_{n_1}(c_1,R),S_{n_2}(c_2,R),\dots,S_{n_t}(c_t,R))$ with $\sum_jn_j=n$. We denote $S_{n_j}(c_j,R)$ by $\Lambda_j$.

For $1\le j\le t$, we define a homomorphism of $R$-modules:
$$\psi_j:M_{n_j}(R)\lra M_n(R),\sum_{1\le u,v\le n_j}x_{uv}e_{uv}\mapsto \sum_{1\le u,v\le n_j}x_{uv}e_{\tau_{j-1}+u,\tau_{j-1}+v},$$
which sends $x=(x_{uv})\in M_{n_j}(R)$ to a $t\times t$ block matrix in which $x$ is in the $(j,j)$-block of size $n_j\times n_j$ and $0$ in $(p,q)$-block of size $n_p\times n_q$ with $(p,q)\ne (j,j)$. Clearly, $(e_{11})\psi_j=e_{\tau_{j-1}+1,\tau_{j-1}+1}\in M_n(R)$, $(a_1)\psi_j(a_2)\psi_j=(a_1a_2)\psi_j$ for $a_1,a_2\in M_{n_j}(R)$ and $(M_{n_j}(R))\psi_j \subseteq \epsilon_j M_n(R)\epsilon_j$, where $1=\sum_{p=1}^{t}\epsilon_p$ is a decomposition of $1$ into pairwise orthogonal idempotents in $M_n(R)$ (see Section \ref{CMA}). Particularly, for $1\le j, l\le t$ and $x\in M_{n_l}(R)$,
$$(*)\quad (e_{11})\psi_j(x)\psi_l=(e_{11})\psi_j\epsilon_j\epsilon_l(x)\psi_l
=\delta_{jl}(e_{11})\psi_l(x)\psi_l=\delta_{jl}(e_{11}x)\psi_l.$$
Since the Jordan blocks in $c_j$ have the same eigenvalues, we have a homomorphism $E_j:M_{n_j}(R)\ra \Lambda_j$ of  $\Lambda_j$-$\Lambda_j$-bimodules as defined in Lemma \ref{stepone} (see the definition of $\bar{E}$). Further, we define a map $$E:M_n(R)\lra S_n(c,R),\;(A_{uv})\mapsto \mbox{diag}((A_{11})E_1,(A_{22})E_2,\dots,(A_{tt})E_t)
=\sum_{p=1}^{t}(A_{pp})E_p\psi_p,$$ where $(A_{uv})$ is a $t\times t$ block matrix with the block $A_{uv}\in M_{n_u\times n_v}(R)$ for $1\leq u,v\leq t$.
Clearly,  $E$ is a homomorphism of $S_n(c,R)$-$S_n(c,R)$-bimodules.

Let $$x_i:=\sum_{p=1}^{t}e_{i,\tau_{p-1}+1},\quad y_i:=\sum_{p=1}^{t}e_{\tau_{p-1}+\lambda_{p1},i}\in M_n(R),\, 1\leq i\leq n.$$

(1) We show that $(E,x_i,y_i)$ is a Frobenius system.
Since $E$ is a homomorphismn of $S_n(c,R)$-$S_n(c,R)$-bimodules, it remains to verify that $\sum_{i=1}^{n} x_i(y_ia)E=a$ and $\sum_{i=1}^{n} (ax_i)Ey_i=a$ for any $a=(a_{pq})\in M_n(R)$. Actually,
$$(**)\quad \sum_{i=1}^{n} x_i(y_ia)E=\sum_{i=1}^{n}\sum_{1\le p, q\le t}e_{i,\tau_{p-1}+1}(e_{\tau_{q-1}+\lambda_{q1},i}  a)E.$$
For $1\le i\le n$ and $1\le p\le t$, $e_{i,\tau_{p-1}+1}(e_{\tau_{q-1}+\lambda_{q1},i}a)E$ is an $n\times n$ matrix with the $i$-th row equal to the $\tau_{p-1}+1$-th row of $(e_{\tau_{q-1}+\lambda_{q1},i}  a)E$ and other rows equal to $0$,
while $e_{\tau_{q-1}+\lambda_{q1},i}a$ is a matrix which has the $(\tau_{q-1}+\lambda_{q1})$-th row equal to the $i$-th row of $a$ and other rows equal to $0$.
Thus $e_{\tau_{q-1}+\lambda_{q1},i}a$ can be written as a $t\times t$ block matrix with $0$ in the $(j,j)$-block of size $n_j\times n_j$ for $1\le j\neq q\le t$, and $\sum_{v=1}^{n_q}a_{i,\tau_{q-1}+v}e_{\lambda_{q1},v}$ in the $(q,q)$-block of size $n_q\times n_q$.
Then, by definition,
$$(e_{\tau_{q-1}+\lambda_{q1},i}a)E=\sum_{j\neq q}(0)E_j\psi_j+(\sum_{v=1}^{n_q}a_{i,\tau_{q-1}+v}e_{\lambda_{q1},v})E_q\psi_q
=(\sum_{v=1}^{n_q}a_{i,\tau_{q-1}+v}e_{\lambda_{q1},v})E_q\psi_q,$$
\begin{align*}
e_{\tau_{p-1}+1,\tau_{p-1}+1}(e_{\tau_{q-1}+\lambda_{q1},i}a)E
&=e_{\tau_{p-1}+1,\tau_{p-1}+1}(\sum_{v=1}^{n_q}a_{i,\tau_{q-1}+v}e_{\lambda_{q1},v})E_q\psi_{q}
\\
&=(e_{11})\psi_p(\sum_{v=1}^{n_q}a_{i,\tau_{q-1}+v}e_{\lambda_{q1},v})E_q\psi_{q} \quad (\mbox{ by } (*))\\
&=\delta_{pq}\,\big(e_{11}(\sum_{v=1}^{n_q}a_{i,\tau_{q-1}+v}e_{\lambda_{q1},v})E_q\big)\psi_{q} \quad (\mbox{ by Lemma }\ref{stepone}(2))\\
&=\delta_{pq}(\sum_{v=1}^{n_q}a_{i,\tau_{q-1}+v}e_{1v})\psi_q=\delta_{pq}\sum_{v=1}^{n_q}a_{i,\tau_{q-1}+v}e_{\tau_{q-1}+1,\tau_{q-1}+v,}
\end{align*}
and therefore the equality $(**)$ runs as follows:
\begin{align*}
\sum_{i=1}^{n} x_i(y_ia)E & = \sum_{i=1}^{n}\sum_{1\le p,q\le t}(e_{i,\tau_{p-1}+1}e_{\tau_{p-1}+1,\tau_{p-1}+1})(e_{\tau_{q-1}+\lambda_{q1},i}a)E\\
&=\sum_{i=1}^{n}\sum_{1\le p,q\le t}e_{i,\tau_{p-1}+1}\big(\delta_{pq}\sum_{v=1}^{n_q}a_{i,\tau_{q-1}+v}e_{\tau_{q-1}+1,\tau_{q-1}+v}\big)\\
&=\sum_{i=1}^{n}\sum_{q=1}^{t}\sum_{v=1}^{n_q}a_{i,\tau_{q-1}+v}e_{i,\tau_{q-1}+v}
\quad (\mbox{ by }[n]=\{\tau_{q-1}+v\in \mathbb{N}\mid 1\le q\le t,1\le v\le n_q\})\\
&=\sum_{1\le i,u\le n}a_{iu}e_{iu} =a.
\end{align*}
Similarly, $\sum_{i=1}^{n} (ax_i)Ey_i=a$. Thus $S_n(c,R)\subseteq M_n(R)$ is a Frobenius extension by Lemma \ref{FROBLEM}(1).

To complete the proof of (1), it remains to prove that the Frobenius extension $S_n(c,R)\subseteq M_n(R)$ is separable. By Lemma \ref{FROBLEM}(2), we have to find an element $d\in M_n(R)$ satisfying the conditions in Lemma \ref{FROBLEM}(2).

Let $b_j:=b_{1j}$ and $\lambda_j:=\lambda_{1j}$ for $1\le j\le s_1$.
Then $n_1=n_{s_1b_{s_1}}$ is the size of $c_1$.
We define $D_i:=\sum_{p=1}^{\lambda_{g(i)}}G_{ii}^p=[1,1,\dots,1]\in M_{\lambda_{g(i)}}(R)$ for $1\leq i\leq m_{s_1},\tilde{d}:=\mbox{diag}(D_1,D_2,\dots,D_{m_{s_1}})=\sum_{i=1}^{m_{s_1}}(D_i)\varphi_{ii}\in M_{n_1}(R)$ and $d:=(\tilde{d})\psi_1\in M_{n}(R)$.
We show that $d$ is a desired element in $M_n(R)$.
In fact, the condition $da=ad$ for all $a\in \Lambda$ is equivalent to saying that $D_iA_{ij}=A_{ij}D_{j}$ holds for all $A_{ij}\in\tilde{\Lambda}_{ij}$ and $1\leq i,j\leq m_s$.
We may write $A_{ij}=\sum_{q=1}^{\theta_{ij}}a_{\theta_{ij}-q+1}G_{ij}^{q}\in\tilde{\Lambda}_{ij}$. It follows from Lemma \ref{multbasis1}(2) and $\lambda_{g(i)}+\theta_{ij}-\lambda_{g(i)}=\theta_{ij}\ge 1$ that
$$D_iA_{ij}= \sum_{u=1}^{\theta_{ij}}\sum_{v=1}^{\theta_{ij}-u+1}a_{(\theta_{ij}-u+1)-v+1}G_{ij}^u
=\sum_{u=1}^{\theta_{ij}}\sum_{w=1}^{\theta_{ij}-u+1}a_wG_{ij}^u=A_{ij}D_{j}$$
for $1\leq i,j\leq m_s$. This means $da=ad$ for all $a\in \Lambda$.

Now we show $\sum_{w=1}^{n}x_wdy_w=I_n$.
It follows from $d=(\tilde{d})\psi_{1} \in \epsilon_1 M_n(R) \epsilon_1$, $e_{w,\tau_{u-1}+1}=e_{w,\tau_{u-1}+1}\epsilon_u$ and $e_{\tau_{u-1}+\lambda_{u1},w}=\epsilon_u e_{\tau_{u-1}+\lambda_{u1},w}$ for $1\le w\le n,1\le u\le t$ that
\begin{align*}
  \sum_{w=1}^{n}x_wdy_w
  &=\sum_{w=1}^{n}(\sum_{u=1}^{t}e_{w,\tau_{u-1}+1})\big((\tilde{d})\psi_{1}\big)
  (\sum_{v=1}^{t}e_{\tau_{v-1}+\lambda_{v1},w})
   =\sum_{w=1}^{n}\sum_{u=1}^{t}e_{w,\tau_{u-1}+1}\epsilon_u (\epsilon_1 (\tilde{d})\psi_{1}\epsilon_1)\sum_{v=1}^{t}\epsilon_v e_{\tau_{v-1}+\lambda_{v1},w}\\
  &=\sum_{w=1}^{n}\sum_{1\le u,v\le t}\delta_{u1}\delta_{1v}e_{w1}(\tilde{d})\psi_{1}e_{\lambda_{11},w}
  =\sum_{w=1}^{n}e_{w1}e_{11}(\tilde{d})\psi_{1}e_{\lambda_{11},\lambda_{11}}e_{\lambda_{11},w.}\\
  &=\sum_{w=1}^{n}e_{w1}(e_{11}\tilde{d}e_{\lambda_{11},\lambda_{11}})\psi_{1}e_{\lambda_{11},w}
  =\sum_{w=1}^{n}e_{w1}(e_{11}\tilde{d}e_{\lambda_{1}\lambda_{1}})\psi_{1}e_{\lambda_{1},w.}
\end{align*}
Further, $I_{n_1}=\sum_{i=1}^{m_{s_1}}f_i$ is a decomposition of $I_{n_1}$ into pairwise orthogonal idempotents in $M_{n_1}(R)$, $(D_i)\varphi_{ii}=f_i(D_i)\varphi_{ii}f_i$, $\tilde{d}=\sum_{i=1}^{m_{s_1}}(D_i)\varphi_{ii}=\sum_{i=1}^{m_{s_1}}f_i(D_i)\varphi_{ii}f_i$ and $e_{11}=e_{11}f_1,e_{\lambda_{1}\lambda_{1}}=f_1e_{\lambda_{1}\lambda_{1}}$ for $1\le i\le m_{s_1}$. Thus
\begin{align*}
  e_{11}\tilde{d}e_{\lambda_{1}\lambda_{1}} & = (e_{11}f_1) (\sum_{i=1}^{m_{s_1}}f_i(D_i)\varphi_{ii}f_i) (f_1e_{\lambda_1\lambda_1})= \sum_{i=1}^{m_{s_1}}\delta_{1i}e_{11}(D_i)\varphi_{ii}e_{\lambda_1\lambda_1}
  =e_{11}(D_1)\varphi_{11}e_{\lambda_1\lambda_1}.
\end{align*}
Since $D_1$ has $1$ in the $(1,\lambda_1)$-entry, $(D_1)\varphi_{11}$ has $1$ in its $(1,\lambda_1)$-entry.
Therefore $e_{11}(D_1)\varphi_{11}e_{\lambda_1\lambda_1}=e_{1\lambda_1}$ and
$$\sum_{w=1}^{n}x_wdy_w=\sum_{w=1}^{n}e_{w1}(e_{11}\tilde{d}e_{\lambda_{1}\lambda_{1}})\psi_{1}e_{\lambda_{1},w}
=\sum_{w=1}^{n}e_{w1}(e_{1\lambda_1})\psi_{1}e_{\lambda_{1},w}=\sum_{w=1}^{n}e_{w1}e_{1\lambda_1}e_{\lambda_{1},w}=\sum_{w=1}^{n}e_{w w}=I_{n.}$$
By Lemma \ref{FROBLEM}(2), the Frobenius extension $\Lambda \subseteq M_n(R)$ is separable.

(2) Now, we prove that $S_n(c,R)\subseteq M_n(R)$ is split if and only if $c={\rm diag}(r_1I_{n_1},r_2I_{n_2},\cdots,r_tI_{n_t})$, where $r_i\in Z(R)$ and $r_i\neq r_j$ for $1\le i\neq j\le t$ and $\sum_{j=1}^tn_j=n.$

If $c={\rm diag}(r_1I_{n_1},r_2I_{n_2},\cdots,r_tI_{n_t})$ with $r_i\in Z(R)$ and $r_i\neq r_j$ for $ i\neq j$, then $S_{n_i}(c_i,R)=M_{n_i}(R)$ and $E_i=id:S_{n_i}(c_i,R)\ra M_{n_i}(R)$ for $1\le i\le t$.
Clearly, $(I_n)E=$ diag $((I_{n_1})E_1,(I_{n_2})E_2,\cdots,(I_{n_t})E_t)=I_n$.
Then $S_n(c,R)\subseteq M_n(R)$ is split by (1) and Lemma \ref{FROBLEM}(2).

Conversely, if $S_n(c,R)\subseteq M_n(R)$ is split, then it follows from (1) and Lemma \ref{FROBLEM}(2) that there exists a $t\times t$ block matrix $z=(Z_{ij})$ with $Z_{ij}\in M_{n_i\times n_j}(R)$ such that $(z)E=I_n$ and $az=za$ for all $a=$ diag $\{a_1,a_2,\cdots,a_t\}\in S_n(c,R)$ with $a_i$ a matrix in $M_{n_i}(R)$ for $1\le i\le t$.
Since $(z)E=$ diag$((Z_{11})E_1,(Z_{22})E_2,\cdots,(Z_{tt})E_t)=I_n$, we have $(Z_{ii})E_{i}=I_{n_i}$ for $1\le i\le t$.
Note that the condition $az=za$ for all $a\in S_n(c,R)$ is equivalent to the condition $a_iZ_{ij}=Z_{ij}a_j$ for all $a_i\in S_{n_i}(c_i,R),a_j\in S_{n_j}(c_j,R),1\le i,j\le t.$
In particular, $Z_{ii}a_i=a_iZ_{ii}$ for all $a_i\in S_{n_i}(c_i,R)$.
By Lemma \ref{stepone}(3), $c_i=r_iI_{n_i}$ for $1\le i\le t$.
By assumption on $c$, the eigenvalue of $c_i$ is not equal to the eigenvalue of $c_j$ for $1\le i\neq j\le t$. Thus $c={\rm diag}(r_1I_{n_1},r_2I_{n_2},\cdots,r_tI_{n_t})$ with $r_i\in Z(R)$ and $r_i\neq r_j$ for $1\le i\neq j\le t$.

(3) If $c={\rm diag}(r_1I_{n_1},r_2I_{n_2},\cdots,r_tI_{n_t})$ and define $c_i=r_iI_{n_i}$, then $S_{n_i}(c_i,R)=M_{n_i}(R)$.
If $R$ is semisimple, then $S_{n_i}(c_i,R)=M_{n_i}(R)$ is semisimple for $1\le i\le t$. Thus $S_n(c,R)$, as a product of these $S_{n_i}(c_i,R)$, is semisimple. Conversely, if $S_n(c,R)$ is semisimple, then $S_{n_i}(c_i,R)$ is semisimple by Lemma \ref{LEMSTEP2}(2) for $1\le i\le t$.
To prove $R$ is semisimple and $c_i=r_iI_{n_i}$ for $1\le i\le t$, it is enough to prove the following claim:

Let $c$ be of the form in $(\dag)$ (see Section 2). If $S_n(c,R)$ is semisimple, then $R$ is semisimple and $c=rI_n$. It follows from Lemma \ref{LEM1STEP1}(1) that $f_iS_n(c,R) f_i\simeq R[x]/(x^{\lambda_{g(i)}})$ is semisimple for $1\leq i\leq m_s$.
This yields that $R$ is semisimple and $\lambda_{g(i)}=1$ for $1\leq i\leq m_s$. In this case, $\lambda_1=1$, $s=1$ and $c=rI_n$. The claim follows.
This also completes the proof of (3).
$\square$

From Theorem \ref{STEPTWO}, we get the corollary.

\begin{Koro}\label{STEPTHR}
Let $k$ be an algebraically closed field.

$(1)$ Every principal centralizer matrix extension over $k$ is a separable Frobenius extension.

$(2)$ If $c\in \Lambda:=M_{n_1}(k)\times M_{n_2}(k)\times\dots\times M_{n_s}(k)$, then $S(c,\Lambda) \subseteq \Lambda$ is a separable Frobenius extension.
\end{Koro}
{\it Proof.} (1) If $k$ is an algebraically closed field, then every square matrix in $M_n(k)$ is a Jordan-similar matrix. Thus Corollary \ref{STEPTHR}(1) follows immediately from Theorem \ref{STEPTWO}.

(2) Let $c=(c_i)\in \Lambda$ with $c_i \in M_{n_i}(R)$ for $1\le i\le n_s$.
Then $S(c,\Lambda) = S(c_1,M_{n_1}(k)) \times S(c_2,M_{n_2}(k)) \times \cdots \times S(c_s,M_{n_s}(k))$.
By (1), $S(c_i,M_{n_i}(k))\subseteq M_{n_i}(k)$ is a Frobenius extension for $1\le i \le s$.
Then $S(c,\Lambda) \subseteq \Lambda$ is a Frobenius extension by Lemma \ref{FROBDP}
$\square$

\medskip
Finally, we remark that, in a general context, the extensions $S_n(C,R)\subseteq M_n(R)$ for $C$ subsets of $M_n(R)$ do not have to be Frobenius extensions.
\begin{Rem} {\rm (1) Let $R$ be a local ring and $n$ be a positive integer such that $nR=0$. If $2\nmid n$, then the extension $S_n(\Sigma_n,R)\subseteq M_n(R)$ is not a Frobenius extension, where $\Sigma_n$ is the symmetric group of degree $n$.

In fact, if $\gamma_n$ denotes the $n\times n$ matrix with all entries equal to $1$, then it follows from $nR=0$ that $\gamma_n^2=0$. Thus $S_n(\Sigma_n,R)= RI_n+R\gamma_n\simeq R[X]/(X^2)$. Since finitely generated projective modules over a local ring must be free and of finite rank, we see that finitely generated nonzero projective $S_n(\Sigma_n,R)$-modules are also free $R$-modules of $R$-rank $2m$ for $m\ge 1$. Due to $2\nmid n$, we deduce that $M_n(R)$ cannot be a projective $S_n(\Sigma_n,R)$-module. Thus $S_n(\Sigma_n,R)\subseteq M_n(R)$ is not a Frobenius extension.

(2) If $C$ contains two matrices in $M_n(R)$ (or if $G$ has no $G$-free point), then $S_n(C,R)\subseteq M_n(R)$ (or $S_n(G,R)\subseteq M_n(R)$) may not be a Frobenius extension.

Indeed, suppose $R$ is a field of characteristic $3$. Due to $\Sigma_3 =\langle (123),(13)\rangle$, it follows for $C:=\{c_{(123)}, c_{(13)}\}\subseteq M_{3}(R)$ that $S_3(C,R)=S_3(\Sigma_3,R)=RI_3+ R\gamma_3$. Then $S_3(C,R)\subseteq M_3(R)$ cannot be a Frobenius extension by (1). Note that $\Sigma_3$ has no free point in $\{1,2,3\}$. \label{rmk}
}\end{Rem}

\section{Cellular algebras\label{CELL}}
In this section, we first recall some basic definition on cellular algebras and then prove Theorem \ref{MAIN2}.

Throughout this note, $R$ stands for a commutative ring with identity. Now, we state the definition of cellular algebras introduced by Graham and Lehrer (see \cite{Graham1996}).
\begin{Def}{\rm \cite{Graham1996}}\label{DEFCELL}
Let $R$ be a commutative ring. A unitary $R$-algebra $A$ is called a \emph{cellular algebra} with cell datum $(P,M,C,\iota)$ if the following conditions are satisfied:

{\rm (C1)} $P$ is a finite partially ordered, and for each $p\in P$, there is associated a finite set $M(p)$ such that $\{C_{ij}^{p}\mid p\in P, i,j\in M(p)\}$ is an $R$-basis of the algebra $A$.

{\rm (C2)} $\iota:A\ra A$ is an $R$-involution (that is an anti-automorphism of $R$-algebra $A$ of order $2$) such that $C_{ij}^{p}$ is mapped to $C_{ji}^{p}$ under $\iota$.

{\rm (C3)} For $a\in A$, $p\in P$, $i,j\in M(p)$, $$aC_{ij}^{p}= \sum_{u\in M(p)}r_{a}(u,i)C_{uj}^{p}+r'$$ where the coefficient $r_{a}(u,i)\in R$ does not depend on $j$ and where $r'$ is a linear combination of basis elements $C^{q}_{st}$ with $q$ strictly smaller than $p$.
\end{Def}

We remark that cellular algebras can also be described in terms of ring-theoretic languages (see \cite{KX} for details). Cellular algebras cover many important examples of algebras such as Hecke algebras, Brauer algebras and Temperley-Lieb algebras, and reduce many problems in representation theory to the ones in linear algebra (see \cite{Graham1996}).

For a cellular algebra $A$ and $p\in P$, we denote by $C^{\le p}$ the $R$-module spanned by all $C^q_{ij}$ with $q\le p$ and $i,j\in M(q)$. By linearization of $P$, we may assume that $P=\{1,2, \cdots, n\}$ with the natural ordering.  Following \cite{KX}, the chain $C^{\le 1}\subset C^{\le 2}\subset\cdots\subset C^{\le n}=A$ ia called a \emph{cell chain} of A.

Note that the cellularity of algebras is not preserved by Morita equivalences. This means that we cannot get cellularity of an algebra by passing to the one of its basic algebra.

Recall that an ideal $I$ of a finite-dimensional algebra $A$ over a field is called a \emph{heredity ideal} if $I=AeA$ for $e^2=e\in A$, $e\rad(A)e=0$ and $_AAeA$ is projective. Following \cite{CPS1988}, a finite-dimensional algebra $A$ is said to be \emph{quasi-hereditary}  if there is a chain of ideals: $0= I_0\subset I_1\subset \cdots \subset I_m=A$ such that $I_i/I_{i-1}$ is a heredity ideal in $A/I_{i-1}$. Such a chain is then called a \emph{heredity chain} of the quasi-hereditary algebra $A$. Quasi-hereditary algebra appears widely in representation theory of Lie algebras and algebraic groups (see \cite{CPS1988}).

\medskip
First, we consider the cellularity of principal centralizer matrix algebras of special form where all blocks have the same eigenvalues.

Let $c={\rm diag}(J_1^{b_1}, J_2^{b_2},\cdots, J_s^{b_s})\in M_{n}(R)$ be a Jordan-block matrix with Jordan blocks $J_{i}=[r,1,0,\cdots,0]\in M_{\lambda_{i}}(R)$ of different sizes $\lambda_i$, $r\in Z(R)$, $1\leq i\leq s$. Here, $J_i^{b_i}$ means that the block $J_i$ appears $b_i$ times. We assume $\lambda_1\textgreater \lambda_2\textgreater \cdots \textgreater \lambda_s$ and keep all notations introduced in Section \ref{CMA}.

For each $p\in [\lambda_1]$, let $l(p)$ be the biggest $l(p) \in [s]$ such that $\lambda_{l(p)}\ge p$. We define $M(p):=[m_{l(p)}]$. Recall that $g(i)$ is the smallest $g(i) \in [s]$ such that $i\le m_{g(i)}$ for $1\le i\le m_s$ and $h(i):=i-m_{g(i)-1}\in [b_{g(i)}]$.

The following lemma is useful in later proofs.
\begin{Lem}\label{llmg}
$(1)$ For $1\le p,q\le \lambda_1$, if $q\le p$, then $l(q)\ge l(p)$ and $M(q)\supseteq M(p)$.

$(2)$ If $1\le u\le s$ and $1\le p\le \lambda_1$, then $p\le \lambda_u$ if and only if $l(p)\ge u$.
In particular, $l(\lambda_u)=u$.

$(3)$ If $1\le i\le m_s$, then $i\in M(p)$ if and only if $g(i)\le l(p)$
if and only if $\lambda_{g(i)}\ge p$.

$(4)$ $g(m_u)=u$ for $1\le u\le s$.
\end{Lem}

{\it Proof.} (1) By definition, we obtain $p\le \lambda_{l(p)}$.
As $q\le p$,  we have $q\le \lambda_{l(p)}$.
The choice of $l(q)$ implies $l(q)\ge l(p)$.
Since $m_1<m_2<\cdots<m_s$, we get $m_{l(q)}\ge m_{l(p)}$.
Therefore $M(q)=[m_{l(q)}]\supseteq [m_{l(p)}]=M(p)$.

(2) If $p\le \lambda_u$, then $l(p)\ge u$ by the choice of $l(p)$.
Conversely, if $l(p)\ge u$, then it follows from $\lambda_1>\lambda_2>\cdots>\lambda_s$ that $\lambda_{l(p)}\le \lambda
_u$.
By the definition of $l(p)$, we have $p\le \lambda_{l(p)}$, and therefore $p\le \lambda_u$.
Specially, if $p=\lambda_u$, then $l(\lambda_u)\ge u$. On the other hand, by the definition of $l(\lambda_u)$, we deduce $\lambda_{l(\lambda_u)}\ge \lambda_u$. It then follows from $\lambda_1>\lambda_2>\cdots>\lambda_s$ that $l(\lambda_u)\le u$. Hence $l(\lambda_u)=u$.

(3) By definition, $g(i)$ is the smallest $g(i) \in [s]$ such that $i\le m_{g(i)}$. If $i\le m_{l(p)}$, that is $i\in M(p)$, then $g(i)\le l(p)$.
Conversely, if $g(i)\le l(p)$, then $i=m_{g(i)-1}+h(i)\le m_{g(i)-1}+b_{g(i)}=m_{g(i)}\le m_{l(p)}$ and $i\in M(p)$.
Hence $i\in M(p)$ if and only if $g(i)\le l(p)$.
By (2),  $p\le\lambda_{g(i)}$ if and only if $l(p)\ge g(i)$.
Thus (3) follows.

(4) By definition, $m_u\le m_{g(m_u)}$.
Since $m_1<m_2<\cdots<m_s$, we conclude $u\le g(m_u)$.
On the other hand, since $g(m_u)$ be the smallest $g(m_u) \in [s]$ such that $m_u\le m_{g(m_u)}$, it follows from $m_u\le m_u$ that $g(m_u)\le u$. Thus $g(m_u)=u$.
$\square$

\medskip
For $i,j\in M(p)$,
we define
$$C_{ij}^p:=\sum_{u=1}^{p}e_{n_{g(i)h(i)}-\lambda_{g(i)}+p-u+1,n_{g(j)h(j)}-u+1}\in M_n(R),$$
that is, $C^p_{ij}= (\sum_{u=1}^{p}e_{p-u+1,\lambda_{g(j)}-u+1})\varphi_{ij}$ is an $m_s\times m_s$ block matrix with $\sum_{u=1}^{p}e_{p-u+1,\lambda_{g(j)}-u+1}$ in the $(i,j)$-block of size $\lambda_{g(i)}\times \lambda_{g(j)}$, and $0$ in the $(k,l)$-block of size $\lambda_{g(k)}\times \lambda_{g(l)}$ if $1\le k,l\le m_s$ and $(k,l)\ne (i,j)$ (see Section \ref{CMA} for the definition of $\varphi_{ij}$).

\begin{Lem} \label{basis} $\Theta :=\{C_{ij}^{p}\mid p\in [\lambda_{1}], i,j\in M(p)\}$ is an $R$-basis of $S_n(c,R)$. \end{Lem}

{\it Proof.} By Lemma \ref{LEM1STEP1}(3), $\Delta:=\{F_{ij}^{p}\mid 1\le i,j\le m_s,1\le p\le \min\{\lambda_{g(i)},\lambda_{g(j)}\}\}$ is an $R$-basis of $S_n(c,R)$. We shall show $\Theta=\Delta$.

Note that for $1\le i,j\le m_{s}$ and $1\le p\le \lambda_1$, if $i,j\in M(p)$ and $p\le \min\{\lambda_{g(i)},\lambda_{g(j)}\}$, then $C_{ij}^p=F_{ij}^p$.
To prove $\Theta=\Delta$, it is enough to show that, for $1\le i,j\le m_{s}$ and $1\le p\le \lambda_1$, we have $i,j\in M(p)$ if and only if $p\le \min\{\lambda_{g(i)},\lambda_{g(j)}\}$ .
But this is clear from Lemma \ref{llmg}(3). Hence
$\Theta=\Delta$ and $\{C_{ij}^{p}\mid p\in P, i,j\in M(p)\}$ is an $R$-basis of $S_n(c,R)$. $\square$

\medskip
For $p,q\in [\lambda_1]$, $i,j\in M(p)$ and $u,v\in M(q)$, by Lemmas \ref{multbasis2} and \ref{basis}, we have the formula
$$(\star) \quad C_{uv}^{q}C_{ij}^{p}=\left\{
	\begin{array}{rcl}
		\delta_{vi}C_{uj}^{p+q-\lambda_{g(v)}} &   & if\;p+q-\lambda_{g(v)}\ge 1,\\
		0\qquad\quad&   & if\;p+q-\lambda_{g(v)}\textless 1.
	\end{array} \right.$$
By using the $R$-basis of $S_n(c,R)$, one  may  define an $R$-linear map
$$\iota: S_n(c,R)\lra S_n(c,R), \; C_{ij}^{p}\mapsto C_{ji}^{p}, \; p\in [\lambda_1], i,j\in M(p).$$
Clearly, $\iota$ is an isomorphism of $R$-modules and $\iota^2=id$.
Given $p,q\in [\lambda_1]$, $i,j\in M(p)$, and $u,v\in M(q)$, if $p+q-\lambda_{g(v)}\ge 1$, then it follows from $(\star)$ that
$$(C_{uv}^{q}C_{ij}^{p})\iota=(\delta_{vi}C_{uj}^{p+q-\lambda_{g(v)}})\iota
 =\delta_{iv}C_{ju}^{q+p-\lambda_{g(i)}}
 =C_{ji}^{p}C_{vu}^{q}
 =(C_{ij}^{p})\iota(C_{uv}^{q})\iota.$$
If $p+q-\lambda_{g(v)}< 1$, then $(C_{uv}^{q}C_{ij}^{p})\iota=(0)\iota=0$ and there always holds $C_{ji}^{p}C_{vu}^{q}=0$ by $(\star)$. This shows $(C_{uv}^{q}C_{ij}^{p})\iota=0=C_{ji}^{p}C_{vu}^{q}$.
In summary, $\iota$ is an anti-automorphism of the $R$-algebra $S_n(c,R)$. Thus $\iota$ is an involution of $S_n(c,R)$.

We remark that the involution $\iota$ is not the transpose of matrices in general.

\begin{Lem}\label{PROPCELL}
Let $c={\rm diag}(J_1^{b_1}, J_2^{b_2},\cdots, J_s^{b_s})\in M_{n}(R)$ be a Jordan-block matrix with Jordan block $J_{i}=[r,1,0,$ $\cdots,0]\in M_{\lambda_{i}}(R)$ for $1\leq i\leq s$ and $r\in Z(R)$,
where $J_i$ appears $b_i$ times
and $\lambda_1\textgreater \lambda_2\textgreater \cdots \textgreater \lambda_s$.
Then
$S_n(c,R)$ is a cellular $R$-algebra with respect to the involution $\iota$.
\end{Lem}
{\it Proof.}
$(1)$ Let $P:=\{1,2, \cdots, \lambda_1\}$ with the natural ordering.
By Lemma \ref{basis}, $\{C_{ij}^{p}\mid p\in P, i,j\in M(p)\}$ with $M(p)=[m_{l(p)}]$ is an $R$-basis of $S_n(c,R)$.

$(2)$ By definition, $(C_{ij}^{p})\iota=(C_{ji}^{p})$ for $p\in P$ and $i,j\in M(p)$.

$(3)$ To verify Definition \ref{DEFCELL}(C3), it is enough to check $(C3)$ for a basis element $a$. Let $p,q\in [\lambda_1]$, $i,j\in M(p)$, and $u,v\in M(q)$. Then
$C_{uv}^{q}C_{ij}^{p}=\delta_{vi}C_{uj}^{p+q-\lambda_{g(v)}}$ if $p+q-\lambda_{g(v)}\ge 1$ and 0 otherwise by $(\star)$.
This means that we have to verify $(C3)$ for the case $p+q-\lambda_{g(v)}\ge 1$. In this case, it follows from $q\le \lambda_{g(v)}$ that the product can be rewritten as
$$C_{uv}^{q}C_{ij}^p= \delta_{v i}\delta_{q, \lambda_{g(v)}}C^p_{u j}+ \delta_{v i}\delta'_{q, \lambda_{g(v)}}C_{uj}^{p+q-\lambda_{g(v)}},$$ where $\delta'_{st}$ stands for the anti-Kronecker symbol, that is, $\delta'_{st}=0$ if $s=t$ and $1$ if $s\ne t$. Note that if $q\ne \lambda_{g(v)}$ then $p+q-\lambda_{g(v)}<p$ and that the coefficients $\delta_{g(v)g(i)}\delta_{vi}\in R$ do not depend on $j$. Hence $C_{uv}^{q}C_{ij}^{p}$ can be expressed in the desired form.
Hence, according to Definition \ref{DEFCELL},  $S_n(c,R)$ is a cellular $R$-algebra with respect to the involution $\iota$. $\square$

In general, we have the result.

\begin{Theo}\label{CELLTHM1} If $R$ is an integral domain and $c$ is an $n\times n$ Jordan-similar matrix over $R$, then $S_n(c,R)$ is a cellular $R$-algebra.
\end{Theo}

{\it Proof.}
Since $c$ is a Jordan-similar matrix in $M_n(R)$, there is a Jordan-block matrix
$d$ such that $c$ is similar to $d$. Thanks to Lemma \ref{PROPCMA}(1), we have $S_n(c,R)\simeq S_n(d,R)$ by an inner automorphism. Let $$d=\mbox{ diag}(J_{11}^{b_{11}},J_{12}^{b_{12}},\cdots,J_{1s_1}^{b_{1s_1}},J_{21}^{b_{21}}, J_{22}^{b_{22}}\cdots,J_{2s_2}^{b_{2s_2}},\cdots,J_{t1}^{b_{t1}},J_{t2}^{b_{t2}}\cdots,J_{ts_t}^{b_{ts_t}})\in M_n(R)$$ with $J_{ij}=[r_i,1,0,\dots,0]\in M_{\lambda_{ij}}(R)$ appearing $b_{ij}$ times for $1\leq j\leq s_i$, $1\leq i\leq t$, $\lambda_{i1}>\lambda_{i2}>\cdots> \lambda_{is_i}$, $r_i\in Z(R)$, and $r_i\neq r_j$ for $i\neq j$.

For $1\leq i\leq t$, we define $n_i:=\sum_{p=1}^{s_i}b_{ip}\lambda_{ip}$, $d_i:=$ diag$(J_{i1},\cdots,J_{i1},\cdots,J_{is_i},\cdots,J_{is_i})\in M_{n_i}(R)$, and $\Lambda_i:=S_{n_i}(d_i,R)$.
Observe that $d_i$ is a Jordan-block matrix with the same eigenvalue $r_i$. It follows from Lemma \ref{LEMSTEP2}(2) that $S_n(d,R)$ is isomorphic to $\Lambda_1\times\Lambda_2\times \dots \times\Lambda_s$ as algebras. Further, by Lemma \ref{PROPCELL}, each $\Lambda_i$ is a cellular $R$-algebra with respect to an involution $\iota_i$.
Now, we define an involution $\iota:=\oplus_{i=1}^t\iota_i:S_n(d,R) \ra S_n(d,R),\; (a_{i}) \mapsto ((a_{i})\iota_i) $ for $a_i\in \Lambda_i, 1\le i\le t.$
Then $S_n(d,R)$ is a cellular $R$-algebra with respect to the involution $\iota$.

It follows from $S_n(c,R) \simeq S_n(d,R)$ that $S_n(c,R)$ is a cellular algebra induced by the cellular structure of $S_n(d,R)$.
$\square$

Consequently, we have the next corollary.

\begin{Koro}\label{CELLTHM3}
Let $R$ be an algebraically closed field.

$(1)$ Every principal centralizer matrix algebra is a cellular $R$-algebra.

$(2)$ If $c\in \Lambda:= M_{n_1}(R)\times M_{n_2}(R)\times\dots\times M_{n_s}(R)$ with $n_i\ge 1$ for all $1\le i\le s$, then $S(c,\Lambda)$ is a cellular algebra.
\end{Koro}
{\it Proof.} $(1)$ Every square matrix over an algebraically closed field is a Jordan-similar matrix. Thus every centralizer matrix algebra is a cellular $R$-algebra by Theorem \ref{CELLTHM1}.

$(2)$ Let $c=(c_1, \cdots, c_s)\in \Lambda$ with $c_i \in M_{n_i}(R)$.
Then $S(c,\Lambda)\simeq S(c_1,M_{n_1}(R)) \times S(c_2,M_{n_2}(R)) \times \cdots \times S(c_s,M_{n_s}(R))$.
Since each $S(c_i,M_{n_i}(R))$ is a cellular algebra by (1), $S(c,\Lambda)$ has a cellular algebra structure induced by $S(c_i,M_{n_i}(R))$ for $1\leq i\leq s$.
$\square$

\medskip
{\bf Proof of Corollary \ref{MAINCOR}:}
(1) is clear from Corollary \ref{STEPTHR}(1).

(2) By a well-known theorem of Maschke, which says that $kG$ is semisimple if and only if the ctaracteristic of $k$ does not divide $|G|$,  we have $kG\simeq M_{n_1}(k)\times M_{n_2}(k)\times\dots\times M_{n_s}(k)$, where $n_1,n_2,\dots,n_s$ are the dimensions of all irreducible representations of $kG$. Now, (2) follows transparently from Corollary \ref{CELLTHM3}(2). $\square$

\medskip
With the help of general theory of cellular algebras, we can parameterize simple modules and describe quasi-heredity of principal centralizer matrix algebras.

We recall the following result on cellular algebras, which is taken from \cite{Graham1996, KX}.

\begin{Lem}\label{lemqhere}
Let $A$ be a cellular algebra over a field $R$ with an involution $\iota$ and cell chain $0=C_0\subset C_1\subset\cdots\subset C_i\subset  C_{\lambda-1}\subset C_{\lambda}=A$. Then the following hold.

$(1)$ There is a natural bijective between isomorphism classes of simple $A$-modules and indices $p\in \{1,2,\cdots,\lambda\}$ such that $C_p^2\not\subset C_{p-1}$.

$(2)$ The given cell chain of $A$ is a heredity chain (making $A$ into a quasi-hereditary algebra) if and only if $C_p^2\not\subset C_{p-1}$ for all $p$ if and only if $\lambda$ equals the number of isomorphism classes of simple modules.
\end{Lem}

The following corollary describes the number of simple modules  and quasi-heredity of principal centralizer matrix algebras.

\begin{Koro} \label{quasi-hered} Let $R$ be a field and $c\in M_n(R)$ be a Jordan-similar matrix of the block type $\{(\lambda_{11}, \lambda_{12},\cdots,\lambda_{1s_1}),$ $\cdots, (\lambda_{t1},\lambda_{t2},\cdots,\lambda_{ts_t})\}.$ Then

$(1)$ $S_n(c,R)$ has $\sum_{i=1}^t{s_i}$ non-isomorphic simple modules.

$(2)$ $S_n(c,R)$ is a quasi-hereditary algebra if and only if $\lambda_{ij}=s_i-j+1$ for $1\le i \le t, 1\le j\le s_i$ if and only if $\lambda_{i1}=s_i$ for $1\le i \le t$.
\end{Koro}

{\it Proof.} First, we prove Corollary \ref{quasi-hered} for the case that $c$ is a Jordan-block matrix with the same eigenvalues, that is, $t=1$ and $c$ has a block type of the form $\{(\lambda_{1},\lambda_{2},\cdots,\lambda_{s})\}$ with $\lambda_1 \textgreater \lambda_2 \textgreater \cdots \textgreater \lambda_s$.
By Lemma \ref{PROPCMA}(1), we can write $c$ = diag$(J_{1}^{b_1}, J_2^{b_2},\cdots,J_{s}^{b_s})\in M_{n}(R)$ as in $(\dag)$.
By Lemma \ref{PROPCELL}, $S_n(c,R)$ is a cellular $R$-algebra with respect to the involution $\iota$. We then have a cell chain $0=C_0\subset C_1\subset\cdots\subset C_p\subset  \cdots\subset C_{\lambda_1}$, where
$C_p$ stands for the $R$-module spanned by all basis elements $C_{ij}^q$ with $q\le p$ and $i,j\in M(q)$.

(1) By Lemma \ref{lemqhere} (1), to prove that $S_n(c,R)$ has exactly $s$ non-isomorphic simple modules, it is sufficient to prove that there are exactly $s$ indices $p\in[\lambda_1]$ such that $C_p^2 \not\subset C_{p-1}$. In the following, we show that $C_p^2 \not\subset C_{p-1}$ if and only if $p=\lambda_{l(p)}$.

We first prove the sufficiency. Suppose $p=\lambda_{l(p)}$. By Lemma \ref{llmg}(4), $g(m_{l(p)})=l(p)$ and $\lambda_{g(m_{l(p)})}=\lambda_{l(p)}=p$. Then it follows from $(\star)$ that $$C_{m_{l(p)},m_{l(p)}}^pC_{m_{l(p)},m_{l(p)}}^p=C_{m_{l(p)},m_{l(p)}}^{p+p-\lambda_{g(m_{l(p)})}}=C_{m_{l(p)},m_{l(p)}}^{p}\in C_p^2\backslash C_{p-1},$$
that is, $C_p^2\not\subset C_{p-1}$.

Now, we show the necessity. Suppose $C_p^2\not\subset C_{p-1}$. According to $(\star)$, there exist $u,v,w\in M(p)$ such that $C_{uv}^{p}C_{vw}^{p}=C_{uw}^{p+p-\lambda_{g(v)}}\not\in C_{p-1}$.
Thus $p+p-\lambda_{g(v)}=p$ and $p=\lambda_{g(v)}$.
By Lemma \ref{llmg}(2), $l(p)=l(\lambda_{g(v)})=g(v)$ and $p=\lambda_{g(v)}=\lambda_{l(p)}$.
Hence, for $p\in [\lambda_1]$, we have $C_p^2 \not\subset C_{p-1}$ if and only if $p=\lambda_{l(p)}$.

By Lemma \ref{llmg}(2), $\lambda_u=\lambda_{l(\lambda_u)}$ for all $1\le u\le s$. Thus
$\{p\in [\lambda_1]\mid p=\lambda_{l(p)} \}=\{\lambda_u\mid 1\le u\le s\}$. Since $\lambda_1>\lambda_2>\cdots>\lambda_s$, the set $\{\lambda_u\mid 1\le u\le s\}$ has exactly $s$ elements.
Thus $S_n(c,R)$ has exactly $s$ non-isomorphic simple $S_n(c,R)$-modules.

(2) By Lemma \ref{lemqhere}(2), that $S_n(c,R)$ is a quasi-hereditary algebra is equivalent to the condition $\lambda_1=s$. Since $\lambda_1>\lambda_2>\cdots>\lambda_s$, we know that $\lambda_1=s$ if and only if $\lambda_i=s-i+1$ for $1\le i\le s$.

Next, we deal with the general case of a Jordan-similar matrix $c\in M_n(R)$. By Lemma \ref{PROPCMA}(1), we may assume $c={\rm diag}(c_1,c_2,\cdots,c_t)$ as in $(\dag\dag)$.
Then $S_n(c,R)$ is isomorphic to $S_{n_1}(c_1,R)\times S_{n_2}(c_2,R)\times\dots\times S_{n_t}(c_t,R))$ as algebras by Lemma \ref{LEMSTEP2}(2).
Since we have shown that $S_{n_i}(c_i,R)$ has $s_i$ non-isomorphic simple $S_{n_i}(c_i,R)$-modules for $1\le i\le t$, the number of non-isomorphic simple $S_n(c,R)$-modules is $\sum_{i=1}^ts_i$.
Clearly, $S_n(c,R)$ is a quasi-hereditary algebra if and only if each $S_{n_i}(c_i,R)$ is a quasi-hereditary algebra for $1\le i \le t$, while $S_{n_i}(c_i,R)$ is a quasi-hereditary algebra if and only if $\lambda_{ij}=s_i-j+1$ for $1\le j\le s_i$ if and only if $\lambda_{i1}=s_i$. This implies that $S_n(c,R)$ is a quasi-hereditary algebra if and only if $\lambda_{ij}=s_i-j+1$ for all $1\le i \le t, 1\le j\le s_i$ if and only if $\lambda_{i1}=s_i$ for $1\le i \le t$. Thus Corollary \ref{quasi-hered} follows. $\square$

\medskip
At this moment, let us display an example to illustrate the results in the paper.

\begin{Bsp} {\rm Let $R$ be a field and $c=$ diag$(J_{1},\dots,J_{2},\dots,J_{s})\in M_{n}(R)$ be a Jordan-block matrix with Jordan blocks $J_{i}=[r,1,0,\cdots,0]\in M_{\lambda_{i}}(R)$, $1\leq i\leq s$,
and $\lambda_1=s, \lambda_2=s-1, \cdots \lambda_{s-1}=2, \lambda_s=1$. In this case, $S_n(c,R)$ is a basic, quasi-hereditary algebra by Corollary \ref{quasi-hered}(2) and Lemma \ref{LEM1STEP1}. }\end{Bsp}

We will work out a presentation of $\Lambda:= S_n(c,R)$ in terms of quiver with relations. By Lemma \ref{LEM1STEP1}, $\Lambda$ has an $R$-basis $F:=\{F_{ij}^{p}\mid 1\le i,j\le s,1\le p\le \mbox{min}\{\lambda
_i,\lambda_j\}\}$ with a complete set $\{F_{ii}^{\lambda_i}\mid 1\le i\le s\}$ of orthogonal primitive idempotent elements, and $\rad(\Lambda)$ has an $R$-basis $F_1:=F \backslash \{F^{\lambda_i}_{ii}\mid 1\le i\le s\}$. Moreover, we prove that $\rad(\Lambda)/\rad^2(\Lambda)$ has an $R$-basis $Q_1:=\{ F^{\lambda_{i}}_{i-1,i}, F^{\lambda_{i}}_{i,i-1} \mid 1< i\le s \}$.

(1) $F_1\backslash Q_1\subseteq \rad^2(\Lambda)$. Suppose $F^{p}_{ij}\in F_1\backslash Q_1$. Then $p\le \mbox{min}\{\lambda_i,\lambda_j\}$. We consider the three cases.

(i) $i < j$. In this case, $\lambda_i > \lambda_j$. If $j= i+1$, then it follows from $F^{\lambda_{i+1}}_{i,i+1}\in Q_1$ that $p < \lambda_{i+1}$. By Lemma \ref{multbasis2}, we have $F_{ij}^p=F_{i,i+1}^p =F_{i,i+1}^{p+1}F_{i+1,i+1}^{\lambda_{i+1}-1}$.
As $F_{i,i+1}^{p+1},F_{i+1,i+1}^{\lambda_{i+1}-1}\in F_1$, we get $F_{ij}^p \in \rad^2(\Lambda)$.
Similarly, if $j>i+1$, then it follows from $F_{i,i+1}^{\lambda_{i+1}}\in F_1$ and $F_{i+1,j}^p\in F_1$ that $F_{ij}^p = F_{i,i+1}^{\lambda_{i+1}}F_{i+1,j}^p\in \rad^2(\Lambda)$.

(ii) $i > j$. In this case, $\lambda_i <\lambda_j$. If $j =i-1$, then it follows from $F^{\lambda_{i}}_{i,i-1} \in Q_1$ that $p < \lambda_{i}$.
By Lemma \ref{multbasis2},  $F^{p}_{i,i-1}=F^{p+1}_{i,i-1}F^{\lambda_{i-1}-1}_{i-1,i-1}\in \rad^2(\Lambda)$.
If $j< i-1$, then $F_{i-1,j}^{\lambda_{i-1}}\in F_1$ and $F^{p}_{ij}=F^{p}_{i,i-1}F^{\lambda_{i-1}}_{i-1,j}\in \rad^2(\Lambda)$.

(iii) $i = j$. In this case, we have $p< \lambda_{i}$ from $F^{\lambda_i}_{ii}\in Q_1$.
Since $\lambda_i=\lambda_{i+1}+1$, this implies $p\le \lambda_{i+1}$.
Then $F_{ii}^p = F_{i,i+1}^pF_{i+1,i}^{\lambda_{i+1}}$ by Lemma \ref{multbasis2}.
It follows from $F_{i,i+1}^p\in F_1$ and $F_{i+1,i}^{\lambda_{i+1}}\in F_1$ that $F_{ij}^p \in \rad^2(\Lambda)$.

(2) No element in $Q_1$ belongs to $\rad^2(\Lambda)$. In fact, for $F_{uv}^{q},F_{ij}^p\in F_1$, the product $F_{uv}^{q}F_{ij}^p$ is either $0$ or again an element of $F_1$.
This implies that $F_{kl}^w\in \rad^2(\Lambda)$ if and only if $F_{kl}^w=F_{uv}^{q}F_{ij}^p$ for some $F_{uv}^{q},F_{ij}^p\in F_1$.
Suppose $F^{\lambda_{k}}_{{k-1},k}=F_{uv}^{q}F_{ij}^p\in Q_1\cap \rad^2(\Lambda)$ for some $F_{uv}^{q},F_{ij}^p\in F_1$.
Then $v=i,u={k-1},j=k, q+p-\lambda_{v}=\lambda_{k}$.
Note that $q\le \mbox{min}\{\lambda_{u},\lambda_{v}\}\le \lambda_{v}$ and $p\le \mbox{min}\{\lambda
_i,\lambda_j\}\le \lambda_{k}$.
Thus $0\le \lambda_{k}-p=q-\lambda_{v}\le 0$.
Therefore $\lambda_{k}=p$, $q=\lambda_{v}$, $\lambda_{k}\le \lambda_{v}$ and $\lambda_{v}\le \lambda_{k-1}$.
As $\lambda_{k}+1=\lambda_{k-1}$, we have $\lambda_{v}=\lambda_{k-1}$ or $\lambda_{v}=\lambda_{k}$.
This means $v={k-1}$ or $v=k$.
If $v={k-1}$, then $F_{uv}^{q}=F_{{k-1}{k-1}}^{\lambda_{k-1}}\not\in F_1$.
If $v=k$, then $i=v=k$  and $F_{ij}^p=F^{\lambda_{k}}_{kk}\not\in F_1$. The both cases contradict to the choices of $F_{uv}^{q}$ and $F_{ij}^p$, respectively.
Thus $F^{\lambda_{k}}_{{k-1},k}\not\in \rad^2(\Lambda)$.
Similarly, $F^{\lambda_{i}}_{i,i-1}\not\in \rad^2(\Lambda)$. Therefore (2) holds and $Q_1$ is an $R$-basis of $\rad(\Lambda)/\rad^2(\Lambda)$.

Now, we define $f_i:=F_{ii}^{\lambda_i}$ for $1\le i\le s$, $\alpha_i:=F_{i,i+1}^{\lambda_{i+1}}$ and $\beta_{i}:=F_{i+1,i}^{\lambda_{i+1}}$ for $1\le i\le s-1$. Then $f_i\alpha_i=\alpha_if_{i+1}$
and $f_{i+1}\beta_i=\beta_i f_i$ for $1\le i < s, \beta_{s-1}\alpha_{s-1}=0 \mbox{ and } \,\alpha_{i}\beta_{i} =\beta_{i-1}\alpha_{i-1}$ for $1< i< s.$
Thus $\Lambda$ is isomorphic to the algebra given by the quiver with relations:

$$\xymatrix{\bullet\ar@<-0.4ex>[r]_{\alpha_{1}}_(0){1}_(1){2}
&\bullet\ar@<-0.4ex>[l]_{\beta_{1}}\ar@<-0.4ex>[r]_{\alpha_{2}}_(1){}
&\bullet\ar@<-0.4ex>[l]_{\beta_{2}}}
\cdots
\xymatrix{\bullet\ar@<-0.4ex>[r]_{\alpha_{s-2}}_(0){}_(1){}
&\bullet\ar@<-0.4ex>[l]_{\beta_{s-2}}\ar@<-0.4ex>[r]_{\alpha_{s-1}}_(1){s}
&\bullet\ar@<-0.4ex>[l]_{\beta_{s-1}}}, \qquad \beta_{s-1}\alpha_{s-1}=0, \; \alpha_{i}\beta_{i} =\beta_{i-1}\alpha_{\i-1}, 1< i< s.
$$
This is actually the Auslander algebra of $R[X]/(X^s)$ which has applications in describing orbits of parabolic subgroups acting on its unipotent radicals (see \cite{Hille}).

To end this section, we propose open questions related to the results in this paper.

\begin{Ques}\label{qu} Let $R$ be an arbitrary unitary ring, and let $\sigma,\tau$ be elements of the symmetric group $\Sigma_n$.

$(1)$ Suppose that $R$ is commutative. Is $S_n(\sigma,R)$ always a cellular $R$-algebra?

$(2)$ Is $S_n(\sigma, R)\subseteq M_n(R)$ always a Frobenius extension?

$(3)$ When are $S_n(\sigma, R)$ and $S_n(\tau, R)$ derived equivalent?

$(4)$ The canonical embedding $R\subseteq M_n(R)$ is a Frobenius extension. How can one parameterize all immediate rings $S$ such that $S\subseteq M_n(R)$ are Frobenius extensions?
\end{Ques}

Note that Question (2) still makes sense, though we have Remark \ref{rmk}. Also, the answer to (3) seems to depend only on numerical properties of the cycle types of $\sigma$ and $\tau$.

\medskip
{\bf Acknowledgement.}
The research work of both authors was partially supported by the National Natural Science Foundation (12031014) and Beijing Natural Science Foundation (1192004). The corresponding author CCX thanks Professor Hourong Qin from Nanjing University for discussing the greatest common divisor matrices which is related to the dimension of $S_n(\sigma,R)$.

\medskip
{\footnotesize

Changchang Xi, School of Mathematical Sciences, Capital Normal University, 100048 Beijing, China; and School of Mathematics and Information Science, Henan Normal University, 453007 Xinxiang, Henan, China

{\tt Email: xicc@cnu.edu.cn}

\medskip
Jinbi Zhang, School of Mathematical Sciences, Capital Normal University, 100048 Beijing, China

{\tt Email: zhangjb@cnu.edu.cn}
}
\end{document}